\documentclass[11pt,draft]{article}
\usepackage{amsmath,amssymb,amsthm}

\title{Free diffusions and Matrix models
 with strictly convex interaction }

\author{Alice Guionnet\thanks{Ecole Normale Sup\'erieure de Lyon,
Unit\'e de Math\'ematiques pures et appliqu\'ees, France
and  Miller institute for
Basic Research in Science, University of California Berkeley. E-mail: aguionne@umpa.ens-lyon.fr.} \and 
D. Shlyakhtenko\thanks{Department of Mathematics, UCLA, 
Los Angeles, CA 90095. E-mail:
shlyakht@math.ucla.edu. Research supported by NSF grants DMS-0355226
and DMS-0555680.} }

\def\half{\frac{1}{2}}
\newtheorem{prop}{Proposition}[section]
\newtheorem{theo}[prop]{Theorem}
\newtheorem{lem}[prop]{Lemma}

\newtheorem{cor}[prop]{Corollary}

\theoremstyle{remark}
\newtheorem{rmk}[prop]{Remark}
\theoremstyle{definition}

\def\P{\mathbb P}
\def\cxmnsa{\C\langle X_1,\ldots,X_m,X_1^*,\ldots,X_m^*\rangle}

\def\C{\mathbb C}
\def\e{\epsilon}
\newcommand{\RR}{\mathbb{R}}

\def\N{{\mathbb N}}

\def\tr{{\mbox{Tr}}}

\def\cxm{{{\mathbb C}\langle X_1,\ldots, X_m\rangle}}

\def\e{\epsilon}

\def\L{\Lambda}

\def\ra{\rightarrow}

\def\Aa{{\cal A}}

\def\Ca{{\cal C}}

\def\La{{\cal L}}
\def\Ma{{\cal M}}

\def\mun{{\hat\mu^N}}

\def\half{\frac{1}{2}}

\def\part{\partial}

\def\ts{\times}

\def\R{\RR}

\begin{document}
\maketitle
\begin{abstract}
We study solutions to the free stochastic differential equation
$dX_t = dS_t - \half DV(X_t)dt$, where $V$ is a locally convex polynomial potential
in $m$ non-commuting variables.  We show that 
for self-adjoint $V$, the law $\mu_V$ of a stationary solution
is the limit law of a random matrix 
model, in which an $m$-tuple of self-adjoint matrices are chosen according to the 
law $\exp(-N \textrm{Tr}(V(A_1,\ldots,A_m)))dA_1\cdots dA_m$.  
We show that if $V=V_\beta$ depends on complex parameters $\beta_1,\ldots,\beta_k$,
then the law $\mu_V$ is analytic in $\beta$ at least for those $\beta$
for which $V_\beta$ is locally convex.  In particular, this gives information
on the region of convergence of the generating function for planar maps.

We show that the solution $dX_t$ has nice convergence properties with respect
to the operator norm.  This allows us to derive several properties of 
$C^*$ and $W^*$ algebras generated by an $m$-tuple with law $\mu_V$.  Among them
is lack of projections, exactness, 
the Haagerup property, and embeddability into the ultrapower of
the hyperfinite II$_1$ factor.  We show that the microstates free entropy
$\chi(\tau_V)$ is finite.

A corollary of these results is the fact that the support of the law of 
any self-adjoint polynomial in $X_1,\ldots,X_n$ under the law $\mu_V$ 
is connected, vastly generalizing the case of a single random matrix.
\end{abstract}
\renewcommand{\refname}{\Large{References}}

\section{Introduction}

There has been a great deal of interest in
 studying matrix integrals in physics since the
work of 't~Hooft who made the connection between 
the problem of enumerating 
maps and estimating integrals of the form
$$Z_N(V)=\int e^{-N\tr(V(X_1,\ldots,X_m))} dX_1\cdots dX_m$$
where $dX$ denotes the 
Lebesgue measure on $N\ts N$ Hermitian matrices
and $\tr$ the non-normalized trace $\tr(A)=\sum_{i=1}^N
A_{ii}$. $V$ is a polynomial in $m$-indeterminates.
 Let us recall that 
a map of genus $g$ is a graph which is embedded
into a surface of genus
$g$ in such a way that the edges do not
intersect and so that disecting the surface along
the edges decomposes it into 
faces, each homeomorphic to a disk. We shall
consider maps with colored edges and enumerate 
them when the degrees of the  vertices, as well
as the ditribution of color of the edges around
each vertex, are  prescribed. The number of
colors will
be $m$, and the colors will be simply refered
by the numbers $\{1,\ldots,m\}$. A vertex with colored half-edges,
an orientation and a distinguished  half-edge,
can be associated bijectivly with a non-commutative monomial 
$q(X)=X_{i_1}\cdots X_{i_p}$ as follows; 
 the  first (or distinguished)  half-edge has  color ${i_1}$,
second has  color ${i_2}$, etc., the
last half-edge having color ${i_p}$.
Such a vertex, equipped with its colored half-edges, distinguished 
half-edge and orientation, will be called a star of
type $q$.  We will
denote by $M_g((q_i,k_i)_{1\le i\le n})$ 
the number of maps with genus $g$ and with
$k_i$ stars of type $q_i$ for $1\le i\le n$, the maps being
constructed by gluing pairwise the half-edges of the stars
(the counting is done up to homomorphism of the surface
and stars are labelled).

't Hooft  
showed that
if $$V(X_1,\ldots, X_m)=W_{(\beta_i,q_i)_{1\le i\le n}}
(X_1,\ldots, X_m)=\frac{1}{2}\sum_{i=1}^m X_i^2+ \sum_{i=1}^n \beta_i
q_i(X_1,\ldots, X_m)$$ where 
$q_i$, $i=1,\ldots,n$ are monomials in $m$-non commutative indeterminates, then
$$\frac{1}{N^2}\log \frac{Z_N(W_{(\beta_i,q_i)_{1\le i\le n}})}{
Z_N(W_{(0,q_i)_{1\le i\le n}}}
=\sum_{g\ge 0}\frac{1}{N^{2g}}
\sum_{k_1,\ldots,k_n\in\N^n} \prod_{1\le i\le n}
\frac{(-\beta_i)^{k_i}}{k_i!} M_g((q_i,k_i)_{1\le i\le n})$$
where the equality holds 
in the sense of formal series.
Differentiating formally this equality, one also finds
that if we consider the Gibbs measure
\begin{equation}\label{eq:RMmeasure}
d\mu^N_{W_{(\beta_i,q_i)_{1\le i\le n}}}
(X_1,\ldots,X_m):=
\frac{1}{Z_N(W_{(\beta_i,q_i)_{1\le i\le n}})}  e^{-N\tr(W_{(\beta_i,q_i)_{1\le i\le n}}(X_1,\ldots,X_m))} dX_1\cdots dX_m
\end{equation}
then we  have, for any monomial $P$,
 the formal expansion

\begin{eqnarray*}
 \bar\mu^N_{W_{(\beta_i,q_i)_{1\le i\le n}}}(P)&:=&\int \frac{1}{N}\tr\left( P(X_1,\ldots,X_m)\right)
d \mu^N_{W_{(\beta_i,q_i)_{1\le i\le n}}}(X_1,\ldots,X_m)\\
&=&\sum_{g\ge 0}\frac{1}{N^{2g}}
\sum_{k_1,\ldots,k_n} \prod_{1\le i\le n}
\frac{(-\beta_i)^{k_i}}{k_i!} M_g((q_i,k_i)_{1\le i\le n}, (P,1)).\\
\end{eqnarray*}
We refer the reader to the survey papers \cite{FGZ, GPW} 
to see  diverse uses of this fact in  the physics literature.

In \cite{APS, EM} for $m=1$ and then 
in \cite{GM,GM2,M} for all $m\in\N$, these formal  equalities 
were  shown to hold in the sense of
large $N$ expansion when $V$ satisfies 
some convexity hypothesis (or one adds a cutoff
to make the integral finite) and the parameters $(\beta_i)_{1\le i\le n}$
are chosen to be small enough. In \cite{GM}, one of the key steps towards this
analysis is to notice that $\bar\mu^N_{W_{(\beta_i,q_i)_{1\le i\le n}}}$ 
converges towards a 
limit, denoted later  $\tau_{W_{(\beta_i,q_i)_{1\le i\le n}}}
$,  which satisfies the so-called Schwinger-Dyson equation 
\begin{equation}\label{SDintro}
\tau\otimes\tau(\partial_i P)=\tau(D_i V  P)
\end{equation}
for all polynomials $P$ and all $i\in\{1,\ldots,m\}$,
and with $V= W_{(\beta_i,q_i)_{1\le i\le n}}$. 
Here, $\partial_i$ and $D_i$ are respectively the 
non-commutative derivative and the cyclic derivative with respect to
the variable $X_i$ (see the next section for a definition). One
then shows that 
 for sufficiently small parameters $\beta_i$, this
equation has a unique solution, which is exactly  
the generating function for planar maps:
$$\tau_{W_{(\beta_i,q_i)_{1\le i\le n}}}(P)=\Ma_{(\beta_i,q_i)_{1\le i\le n}}
(P):=\sum_{k_1,\ldots,k_n} \prod_{1\le i\le n}
\frac{(-\beta_i)^{k_i}}{k_i!} M_0((q_i,k_i)_{1\le i\le n}, (P,1)).$$

It is natural 
to wonder how
to use these expansions to study the numbers $M_g((q_i, k_i)_{1\le i\le n} )$
and in particular their asymptotics as $k_1,\ldots,k_n$
go to infinity.  The answer to this
question is still open in such a general
context. It is, however, quite well understood
in the case $m=1$.
Let us highlight this point in the case of quadrangulations,
corresponding to the potential $V(x)=\beta x^4$, $\beta\in\R$, 
even though 
the enumeration of quadrangulations
was achieved by direct combinatorial arguments
by Tutte \cite{Tu} long ago.
In this case, the Schwinger-Dyson equation \eqref{SDintro}, taken at $P(x)=(z-x)^{-1}$
shows that the Cauchy transform
$$G(z)=\int \frac{1}{z-x} 
d\tau_{W_{\beta,x^4}}(x)=\int \frac{1}{z-x} d
\Ma_{\beta,x^4}
(x)$$
satisfies an algebraic equation of degree two
$$G(z)^2=4\beta z^3 G(z)+zG(z)+P(z)\quad
\mbox{ with }
P(z)=4\beta \int\frac{x^3-z^3}{x-z} 
d\tau_{W_{\beta,x^4}}(x)-1.$$
We can solve this equation in terms of $P(z)$
which is a polynomial of degree two with two unknown coefficients;
we then find that $G$ is given by a polynomial plus the square root 
of a polynomial of degree six, that we denote $Q$, with two unknown coefficients.
Until this point, all the arguments follow 
the induction relations already found by Tutte.
However, the difference now is that
we know that $\Ma_{\beta,x^4}
=\tau_{W_{\beta,x^4}}$ 
is a probability measure 
on $\R$. Assume we can argue that for sufficiently small $\beta$,
the support of $\tau_{W_{\beta,x^4}}$ is connected. Then, we see that,
because this means that $G$ is analytic outside an
interval, the polynomial $Q$ must have two double roots.
This actually determines $Q$, and thus $G$, uniquely.
Since $G$ is the generating function
for planar maps, we are done. 
Hence, we see in this context that the {\em a priori}
information that $\Ma_{(\beta_i,q_i)_{1\le i\le n}}
(P)$ is the Cauchy transform of 
a measure on the real line (which is not clear from its definition
as a generating function of maps), and  with
a connected support,
is enough to conclude.

The goal of this article is to push forward
the analysis of the limiting tracial state 
$\tau_{W_{(\beta_i,q_i)_{1\le i\le n}}}$
in the  multi-matrix  context. We prove in particular that
when ${W_{(\beta_i,q_i)_{1\le i\le n}}}$ satisfies a certain local convexity
property (see \eqref{localconvex}), 
the support of the limiting spectral measure of the  random matrices
with law  $\mu^N_{W_{(\beta_i,q_i)_{1\le i\le n}}}$ is connected.
In fact, the same is true for an arbitrary non-commutative polynomial
in the random 
matrices.  
Note that this information
is enough to solve the enumeration
problem when $m=1$ as we have seen above
for quadrangulations  (though connectivity of support
can in this case 
 be proved by other techniques, see e.g
\cite{deift}).

The tracial states  $\tau=\tau_V$ under consideration
will be solution to the
Schwinger-Dyson equation \eqref{SDintro}
for some general potential $V$.
Non-commutative laws arising as limits of laws of random matrix
models given by \eqref{eq:RMmeasure} have also 
naturally appeared in free probability 
theory.  There, the fact
that  they satisfy 
the  Schwinger-Dyson type equation is restated
as the fact that  the 
free conjugate variables of the law are equal to the cyclic gradient
of a polynomial potential, see also \cite{CDG2,biane}.

In the  multi-matrix setting,  uniqueness of the solution 
to the Schwinger-Dyson equation
is unclear in general. It was proved in \cite{GM}
that, when $V=W_{(\beta_i,q_i)_{1\le i\le n}}$,
 there exists a unique solution 
such that $|\tau(X_{i_1}\cdots X_{i_k})|\le R^k$
for all choices of $i_j\in\{1,\cdots,m\}$ and all $k=1,2,\ldots$, provided the
$\beta_i$'s are sufficiently small. 
In this paper we define a notion  of locally strictly  convex potential,
which generalizes to non-commutative variables the standard notion
of local convexity for functions on the
real line. One of the central result
of this paper will be the uniqueness of
the solution to  Schwinger-Dyson equation for
locally  strictly convex potentials $V$, when 
the domain of strict convexity is large enough.

We also show that if $V=V_\beta= \sum \beta_i q_i$
with some monomials $q_i$ and $\beta=(\beta_i)_{1\le i\le n}
$, $\beta\ra \tau_{V_\beta}(P)$
is analytic in the whole region of local convexity 
of $V_\beta$, for any monomial $P$. Because 
$$ \Ma_{{(\beta_i,q_i)_{1\le i\le n}}}(P)=\tau_{W_{(q_i,\beta_i)_{1\le i\le
n}}} =\sum_{k_1,\ldots,k_n} \prod_{1\le i\le n}
\frac{(-\beta_i)^{k_i}}{k_i!} M_0((q_i,k_i)_{1\le i\le n}, (P,1)),$$
this result shows that there is no breaking of analyticity 
of $(\beta_i)_{1\le i\le n}\ra  \Ma_{{(\beta_i,q_i)_{1\le i\le n}}}(P)$
 in
the domain where $ W_{(q_i,\beta_i)_{1\le i\le
n}}$ stays locally strictly convex, and thus
provides valuable information on the
asymptotics of the numbers $ M_0((q_i,k_i)_{1\le i\le n}, (P,1))$.

Of particular interest are the $C^*$-algebras $(A_V,\tau_V)$
 and the $W^*$-algebras $(M_V,\tau_V)$ generated by 
operators having a law satisfying the Schwinger-Dyson equation with 
a fixed locally convex potential $V$.  
We derive several properties of such algebras, showing that they are 
projectionless, exact \cite{wassermann} and possess the compact approximation property of 
Haagerup \cite{haagerup:compact}, and can be embedded into the ultrapower of the hyperfinite II$_1$
factor.  We also show that the algebras $M_V$ are factors and that the generating 
operators have a law with finite free entropy.  This has as consequences a number of properties
of the algebras $M_V$ (including primeness and lack of Cartan subalgebras), 
see \cite{dvv:entropy3} and \cite{ge:prime}.
  All of these properties are similar to (and are often derived from)
the corresponding properties of free group factors (which are von Neumann
algebras associated to a quadratic potential $V(X_1,\ldots,X_n)=\half\sum_i X_i^2$).
This adds evidence towards a positive answer to Voiculescu's question 
of whether $M_V$ is isomorphic to a free group factor for a fairly 
arbitrary potential $V$.

While somewhat technical, it should be noted that studying properties of
$A_V$ and $M_V$ is of substantial interest.  
Indeed, connectivity of support of limit
distributions of random matrix models is directly related to the lack of 
projections in the $C^*$-algebra $A_V$ (which in turn is derived from the 
famous result of Pimsner and Voiculescu \cite{PV}, essentially
 dealing with the quadratic
$V$).

Finally, in the remaining sections of the paper we show that the random 
matrices following the law \eqref{eq:RMmeasure} give a very good approximation
to the non-commutative law $\tau_V$.  If $V$ is locally convex, we show that
the $\limsup$ in the definition of the microstates free entropy of $\tau_V$ 
\cite{dvv:entropy2} can be replaced by a limit.  In the case that $V$
is (globally) convex, we show that operator norms of arbitrary polynomials
in such random matrices almost surely approximate the operator norms
of such polynomials computed in the $C^*$-algebra $A_V$.  This extends
the results of \cite{HT}, which correspond to the case of quadratic $V$
(our proof, though, relies on their result).

The main technical tools used in the present paper involve the study of 
free Langevin-type diffusion
and its convergence to
a stationary measure
which  corresponds to the limit law  $\tau_V$ 
of random matrices following the 
measure \eqref{eq:RMmeasure}.  We extend some of the results of \cite{BS2}
to the setting of locally convex potentials (see below).  This way, we 
are able to show that operators having a specific limit law can be 
approximated in the operator norm by continuous functions of free 
Brownian motion.  This enables us to carry over a number of properties
of the algebra generated by free Brownian motion to the algebras $A_V$
and $M_V$.

\subsection{Definitions, notations and statment of the results.}
Let us now state more precisely our setup and results.

We let $\cxm$ be the set of
polynomials in $m$ non-commutative 
variables $(X_1,\ldots,X_m)$. 
We shall not assume in general 
that $(X_1,\ldots,X_m)$ are self-adjoint
but let $(X_1^*,\ldots,X_m^*)$
be their adjoints for some involution $*$.
We denote $\cxmnsa$ the set of polynomials in 
the non-commutative variables
$(X_1,\ldots,X_m,X_1^*,\ldots, X_m^*)$.
This set is endowed with the linear involution
so that
$$(X_{i_1}^{\varepsilon_1} \cdots X_{i_k}^{\varepsilon_k})^*=
X_{i_k}^{-\varepsilon_k}
X_{i_{k-1}}^{-\varepsilon_{k-1}}\cdots X_{i_1}^{-\varepsilon_{1}}$$
where we denoted in short $X_i^1=X_i$ and $X_i^{-1}=X_i^*$
and the $(\varepsilon_1,\ldots,\varepsilon_k)$
belong to $\{-1,1\}^k$. We shall denote below
for two sets of non-commutative
variables $(X_1,\ldots,X_m)$ and
$(Y_1,\ldots,Y_m)$ and  an involution
$*$ 
$$X.Y=\frac{1}{2}\sum_{i=1}^m (X_iY_i^*+Y_iX_i^*).$$
$\|\cdot\|_\infty$ will denote an operator norm
such that the completion of $(\cxmnsa,*)$ 
for this norm is a $C^*$-algebra. For a $m$-dimensional vector $X=(X_1,\ldots,
X_m)$ we denote in short $\|X\|_\infty=\max_{1\le i\le m}
\|X_i\|_\infty$.

We let $V\in\cxm$ be a polynomial in $m$ non-commutative 
variables.
 We will say that  
$V$  is $(c,M)$ convex iff 
for any  $m$-tuples of non-commutative  variables 
 $X=(X_1,\dots,X_m)$ and $Y=(Y_1,\dots,Y_m)$ 
in some $C^*$-algebra $(\cal A,\|\cdot\|_\infty)$
satisfying $\Vert X_i\Vert_\infty, \Vert Y_i\Vert_\infty \leq M$, $i=1,\ldots,m$
we have
\begin{equation}\label{localconvex}
[DV(X)-DV(Y)].(X-Y)\ge c (X-Y).(X-Y)
\end{equation}
where the inequality is understood in the sense of operators
($X\ge Y$ iff $X-Y$ is self adjoint
and has non negative spectrum).
 $D=(D_1,\ldots,D_m)$
denotes the cyclic gradient which is
linear and given, for any monomial $P$,
by
$$D_iP=\sum_{P=P_1X_i P_2} P_2P_1.$$
Later, we shall also need the non-commutative gradient $\partial$
which is given, for any monomial $P$,
any $i\in \{1,\ldots,m\}$
by
$$\partial_iP=\sum_{P=P_1X_i P_2} P_1\otimes P_2.$$
We occasionally shall consider
polynomials in $\cxmnsa$; in that case we
 extend $\partial_i$ and $D_i$ by setting 
$\partial_i X_j=1_{i=j}1\otimes 1$ and $\partial_i X_j^*=0\otimes 0$
whereas we have also the derivative $\partial_{i,*}$ and $D_{i,*}$
with respect to $X_i^*$ which satisfy $\partial_{i,*}X_j=0\otimes 0$
but $\partial_{i,*} X_j^*=1\otimes 1$,  and extending
by linearity and Leibnitz rule.

In the case that $X_1,\ldots,X_n$ are self-adjoint, we shall make the 
convention that $\partial_i X_j = 1_{i=j} 1\otimes 1$, while
$\partial_{i,*}X_j = 0$ for all $j$ (in other words, we shall continue
to think of all quantities as functions of $X_1,\ldots,X_n$ alone).

Note that by taking $X=Y+\e Z$ and
letting $\e$ going to zero
we find that for any bounded operator $Z$
and any operator $Y$ with norm strictly smaller than 
$M$,
the condition that $V$ is $(c,M)$-convex implies that
$$\sum_{i=1}^m \sum_{j=1}^m  \left(
\partial_i D_j V(Y)\sharp Z_i \times Z_j^*+
Z_j\times (\partial_i D_j V(Y)\sharp Z_i)^*\right)\ge 2cZ.Z.$$

We shall say that $V\in\cxm$  is self-adjoint  iff
for any self-adjoint  variables 
 $X=(X_1,\dots,X_m)$, $[V(X_1,\ldots,X_m)]^*= V(X_1,\ldots,X_m)$.
$V$ is self-adjoint $(c,M)$-convex
if $V$ is self-adjoint 
and the above holds once restricted to
self-adjoint variables, i.e 
for any  $m$-tuples of self-adjoint variables 
 $X=(X_1,\dots,X_m)$ and $Y=(Y_1,\dots,Y_m)$ 
living in some $C^*$-algebra $(\cal A,\|\cdot\|_\infty)$
which are bounded in norm by $M$,
we have
$$[DV(X)-DV(Y)].(X-Y) \ge c (X-Y).(X-Y)
$$
In this case  $X.Y=\{X,Y\}:=\frac{1}{2}\sum_{i=1}^m (X_i Y_i +Y_i X_i)$
is simply the  anti-commutator of $X$ and $Y$.

If we specialize this assumption to
matrices and consider
$\Aa$ to be the algebra 
of $N\ts N$ matrices
with complex entries
equipped with the usual involution $(A^*)_{ij}= \bar A_{ji}$
and the spectral norm $\|\cdot\|_\infty$,
we find that if $V$ is self-adjoint $(c,M)$-convex, 
$(X_{ij})_{i\le j}\ra \tr [ V(X)]$ is strictly convex 
on the set of entries where $X$ is Hermitian 
and with spectral radius bounded by $M$ since 
\begin{multline*}
\tr V(X)-\tr V(Y)=\int_0^1 \tr (DV(\alpha X+(1-\alpha)Y).(X-Y)) 
d\alpha\\
 =  \tr (DV(Y).(X-Y)) 
+\int_0^1 \tr ([DV(\alpha X+(1-\alpha)Y)-DV(X)].(X-Y)) 
d\alpha\\
\ge\tr (DV(Y).(X-Y)) +\frac{c}{2} \tr((X-Y)^2).
\end{multline*}
Taking $Y=(X+Z)/2$ and $X$ to be $X$ or $Z$
and summing the resulting inequalities
gives 
$$\tr V(X)+\tr V(Z)-2\tr V\left(\frac{X+Z}{2}\right)
\ge c \tr((X-Z)^2)$$
and 
hence the Hessian of $(X_{ij})_{i\le j}\ra \tr V(X)$
is bounded below by $c I$, at least on matrices $X$
with norm bounded by $M$.

This kind of hypothesis was shown to be very useful
in \cite{GM}.
The interest
in relaxing the
hypothesis of convexity to
hold in a bounded domain is related with 
matrix models where $V=\frac{1}{2}X.X+W$
with $W=\sum_{i=1}^n  \beta_i q_i$
for some monomials and complex parameters $(\beta_i)_{1\le i\le n}\in\C^n$.
It is clear now that for all $M$ finite
we can choose the $\beta_i$'s
sufficiently small so that
$V$ is $(1/2,M)$-convex (whereas it 
would not work with no bounds).
Indeed, in that case
$$(DV(X)-DV(Y)).(X-Y)=(X-Y).(X-Y)+(DW(X)-DW(Y)).(X-Y)$$
But when the norms of $X$ and $Y$ are bounded
by $M$,
$$|(DW(X)-DW(Y)).(X-Y)|\le C(M)\max_{i}|\beta_i| (X-Y).(X-Y)$$
with $C(M)$ a constant which only
depends on $M$ and the $q_i$.
Hence, we can now choose $t$ small enough
so that $C(M)\max_{i}|t_i|<1/2$ and so $V$ 
is then $(1/2,M)$-convex.
 This is analogous 
to what
was done in \cite{GM} 
in case of non-convex interaction; it was 
shown that then if one  adds a cut-off, the
large $N$ expansion is still valid
provided the parameters in $W$ 
are small enough.

Hereafter we assume that $V$ is $(c,M)$-convex.
We let $(\cal A,*, \phi)$
be a non-commutative probability space 
generated by a free Brownian
motion $S$  (we refer to \cite{BS2,BS} for 
an introduction
to free Brownian motion and its related free It\^o calculus).
We shall denote by $\|\cdot\|_\infty$
the operator norm in  $(\cal A,\phi)$.

We prove (see Lemma \ref{exist} and Theorem \ref{existtheo}):
\begin{theo}\label{main}
Let $V$ be a $(c,M)$-convex polynomial in $X_1,\ldots,X_m$.

Then there exist 
$M_0=M_0(c, \|DV(0).DV(0)\|_\infty)$, 
$B_0=B_0(c, \|DV(0).DV(0)\|_\infty)$, and $b=b (c,\|DV(0).DV(0)\|_\infty, M) \ge B_0$
finite constants, 
so that whenever $M\geq M_0$  and $Z$ is an $m$-tuple with $\Vert Z \Vert < b$,
\renewcommand{\theenumi}{\roman{enumi}}
\begin{enumerate}
\item 
 There exists a unique solution  $X_t^Z$
to
\begin{equation}\label{sde0}
dX^Z_t=dS_t-\half DV(X_t^Z)dt, \qquad t\in [0,+\infty),
\end{equation} with the initial data $X_0^Z=Z$.
Moreover, in this case,
\begin{eqnarray*}
\Vert X_t^Z\Vert_\infty \leq M,&\qquad& \forall\, t\in [0,+\infty),\\
\limsup_{t\to\infty} \Vert X_t^Z\Vert_\infty\leq B_0,&& \\
X_t^Z\in C^*(Z, S_q : q\in [0,t]),&\qquad&\forall\, t\in[0,+\infty). 
\end{eqnarray*} 
\item $\Vert X_t^Z - X_t^0 \Vert_\infty \to 0\qquad\textrm{as $t\to\infty$}.$
\item The law of $((X_t^Z)^*,X_t^Z)$
converges to a stationary 
law  $\mu_V\in\cxmnsa'$ as $t$ goes to infinity.
$\mu_V$ is the non-commutative  law of $m$ 
variables uniformly bounded by $B_0$.
If $V$ is self-adjoint, $\mu_V$ is the law of $m$ self-adjointvariables. 
\item
The restriction $\tau_V= \mu_V|_\cxm$ 
of $\mu_V$ to $\cxm$ satisfies the
Schwinger-Dyson
equation which states that for all $P\in\cxm$ and all $i\in\{1,\ldots,m\}$
\begin{equation}\label{SD00}
\tau_V\otimes\tau_V(  \partial_i  P)=\tau_V(D_i V P).
\end{equation}
Moreover, if $\nu$ is the law of $m$ variables whose uniform norm is bounded by $b$ and
$\nu|_\cxm$ satisfies \eqref{SD00},
then $\nu|_\cxm=\tau_V$.
\end{enumerate}
\end{theo}
The Schwinger-Dyson equation \eqref{SD00}
is exactly the same as the one 
which characterized the enumeration
of maps \eqref{SDintro} from which we deduce 
that $\Ma_{{(\beta_i,q_i)_{1\le i\le n}}}=\tau_{\frac{1}{2}\sum_{i=1}^m X_i^2
 +\sum \beta_i q_i}$ at least for sufficiently small
polynomials $V$. This allows to 
give some information on the domain of
analyticity of ${\bf \beta}=(\beta_i)_{1\le i\le n}
\ra\Ma_{{(\beta_i,q_i)_{1\le i\le n}}}(P)=
 \Ma_{\bf \beta}(P)$ (see section \ref{analytic}).
\begin{theo}
Let $V=V_{\beta}=\sum_{i=1}^n \beta_i q_i$ be a polynomial, where
${\bf \beta}=(\beta_i)_{1\le i\le n}$ are (complex) parameters 
and $(q_i)_{1\le i\le n}$ are monomials. For $c,M$ positive real
numbers and $M\ge M_0$ with $M_0$ as in Theorem  \ref{main}, 
let $T(c,M
)\subset \C^n$ 
be the interior 
of the subset of parameters ${\bf \beta}=(\beta_i)_{1\le i\le n}$
such that $V_\beta$ is $(c,M)$-convex.

Then , for any $P\in\cxm$,
$\beta\in T(c,M
) \ra \tau_{V_\beta}(P)$ is analytic.
In particular, $\beta\ra \Ma_\beta(P)$
extends analytically to the interior of the set of $\beta_i$'s
where $\frac{1}{2}\sum_{i=1}^m X_i^2 +\sum \beta_i q_i$
is $(c,M)$-convex for $M\ge M_0(c)$ .
\end{theo}
The laws $\mu_V$ are interesting in their own. In the free 
probability language we have proved
the following.

\begin{theo}
Let $V$ be a  $(c,M)$-convex potential 
with $M\ge M_0$ the constant of Theorem \ref{main}.
If $V$ is self-adjoint, there exists a  non-commutative law
in $\cxmnsa'$
with conjuguate variable $(D_i V)_{1\le i\le m}$.
There exists at most one such law satisfying the
additional constraint to be
the law of variables bounded by $b$.
For non self-adjoint
potential, there  a unique law  $\mu_V\in\cxmnsa'$
which satisfies for all $P\in\cxmnsa$ and all $i\in\{1,\ldots,m\}$

\begin{equation}\label{SDtot}
\sum_{i=1}^m \mu_V\otimes\mu_V\left(  (\partial_i+\partial_{i,*})
(D_i+D_{i,*}) P\right)=\sum_{i=1}^m\mu_V\left( D_i V D_iP+ (D_i V)^*
D_{i,*} P\right)
\end{equation}
with $(D_{i,*},\partial_{i,*})$ the non-commutative deriavtives
with respect to $X_i^*$. There exists at most one such law 
satisfying the
additional constraint to be
the law of variables bounded by $b$.

\end{theo}
This statement is deduced from \eqref{existtheo}(2).
Moreover, we let $Z$ be an $m$-tuple of operators with
law $\mu_V$. We shall prove that the $C^*$-algebra 
and the von Neumann algebra generated
by $Z$ have  many properties in common
with the one generated by a semi-circular system.

\begin{theo} Assume that $V$ is $(c,M)$-convex with 
$M\ge M_0$ the constant of Theorem \ref{main}. The
 $C^*$-algebra generated by $Z$
is exact, projectionless and its  associated
von Neumann algebra has the Haagerup approximation property and admits
and embedding into the ultrapower of the hyperfinite II$_1$ factor. 
\end{theo}
In particular we have
\begin{cor} Assume that $V$ is $(c,M)$-convex with 
$M\ge M_0$ the constant of Theorem \ref{main}. 
Let $Z$ be any $m$-tuple of $b$-bounded variables
with the unique law $\mu_V$ satisfying
\eqref{SDintro}.
(a) The algebra $C^*(Z)$ has no non-trivial projections.  
(b) The spectrum of any non-commutative
*-polynomial $P$ in the $m$-tuple $Z$ is connected (in the case that 
$P(Z)$ is normal, this means that the support of its spectral measure is 
connected).  
(c) If $P$ is any polynomial in $Z$ whose value is self-adjoint, then 
the probability measure given by the law of $P(Z)$ has connected support.
\end{cor}

\section{Existence of free diffusions and convergence
to their stationary measure}

We shall show  that if $M$ is chosen
large enough (depending on $c$ and
$ \|DV(0).DV(0)\|_\infty$),
we can build a bounded solution
to some free stochastic
differential equation with drift $DV$ 
provided that $V$ is $(c,M)$-convex. 
This generalizes Langevin dynamics to the
context of operators.
\begin{lem}\label{exist}
Let $V$ be a $(c,M)$-convex polynomial in $X_1,\ldots,X_m$.
Then there exist finite 
constants 
\begin{gather*}
M_0=M_0(c, \|DV(0).DV(0)\|_\infty),\qquad 
B_0=B_0(c, \|DV(0).DV(0)\|_\infty)\\
b=b (c,\|DV(0).DV(0)\|_\infty, M) \ge B_0
\end{gather*}
so that if $M\geq M_0$,
and $Z$ is any $m$-tuple with $\Vert Z \Vert < b$,
 there exists a unique solution  $X_t$
to
\begin{equation}\label{sde}
dX_t=dS_t-\half DV(X_t)dt, \qquad t\in [0,+\infty).
\end{equation} with the initial data $X_0=Z$.
Moreover, in this case,
\begin{eqnarray*}
\Vert X_t\Vert \leq M,&\qquad& \forall\, t\in [0,+\infty),\\
\limsup_{t\to\infty} \Vert X_t\Vert_\infty\leq B_0,&& \\
X_t\in C^*(Z, S_q : q\in [0,+\infty)),&\qquad&\forall\, t\in[0,+\infty). 
\end{eqnarray*} 

If $V$ is self-adjoint  $(c,M)$-convex,
the above results hold under the additional
assumption that $X_0=Z$ is self-adjoint.
In this case, $X_t$ remains self-adjoint for all $t\ge 0$.
\end{lem}
{\bf Proof.}
We remind the reader (cf. \cite{BS2}) that 
if $DV$ is uniformly Lipschitz for
the uniform norm, 
the existence and uniqueness 
to \eqref{sde} is clear by the following 
Picard argument. For the existence we consider
the sequence $X_t^n$, $n=0,1,2,\ldots$ constructed recursively as follows.
Set $X^0_t=Z$ for
all $t$, and having defined $X^n_t$, 
let $X^{n+1}_t$ be given by the
equation
$$dX_t^{
n+1}=dS_t-\half DV(X_t^n)dt$$
with the initial condition $X^{n+1}_0=Z$.

Subtracting the equations for $X_t^{n+1}$ and $X_t^n$ from each other, 
we get for all $t\ge 0$, 
$$\|X_t^{n+1}-X^n_t\|_\infty\le\half
\int_0^t \|DV(X^n_s)-DV(X^{n-1}_s)\|_\infty
ds\le \half \|DV\|_\La \int_0^t \|X^n_s-X^{n-1}_s\|_\infty
ds$$
where $\|DV\|_\La$ denotes the Lipschitz norm of $DV$
$$\|DV\|_\La = \inf\{C: \Vert DV(Y)-DV(Y')\Vert_\infty 
\leq C \Vert Y-Y'\Vert_\infty\}.$$
Iterating, we deduce that 
$$\|X_t^{n+1}-X^n_t\|_\infty\le  \left(\frac{\|DV\|_\La}{2}\right)^n \frac{t^n}{n!}\sup_{u\le t}
\|S_u+u DV(Z)\|_\infty\le  \left(\frac{\|DV\|_\La}{2}\right)^n \frac{t^n}{n!}(2\sqrt{t}
+\frac{t}{2} \|DV(Z)\|_\infty )$$
which proves norm convergence of $X_t^n$ (note that $X_0^n = Z$ for all
$n$).  Moreover, this limit satisfies 
\eqref{sde}.   We also see that $X_t^n \in C^*(Z, S_q : q\geq 0)$ and 
therefore the limit $X_t \in C^*(Z, S_q: q\geq 0)$ as well.

The proof of uniqueness follows the same lines since 
any two 
solutions $X_t,Y_t$  satisfy
$$\|X_t-Y_t\|_\infty\le \half \|DV\|_\La \int_0^t \|X_s-Y_s\|_\infty
ds$$
which proves that $X=Y$ by Gronwall's argument.

If $V$ is a  self-adjoint polynomial,
then $D_iV$ is also self-adjoint for
all $i\in\{1,\ldots,m\}$.
Therefore, since $S_t$ is self-adjoint for all
$t\ge 0$, we deduce by induction that $X^n_t$
is self-adjoint for all $n\ge 0$ and
all $t\ge 0$ and
so its limit $X_t$ is
also self-adjoint. 
Hence the solution  $X_t$
of the free stochastic equation  is self-adjoint for
all times. 

We now return to the case of a $(c,M)$ convex $V$.
By Lemma 3.2 in \cite{BS2}, we can construct a new function
$$f_R(X) = \half DV(X) h(\sum_i \Vert X_i\Vert_\infty)$$ so that $f_R$ 
is uniformly Lipschitz
and  $f_R(X)=DV(X)$ if $\sum \Vert X_j \Vert_\infty \leq R$.
The Picard argument above implies existence and uniqueness of a solution
$X^R_t$ 
to $$dX^R_t = dS_t + f_R(X^R_t) dt$$ for all times $t$.  Clearly, if we show
that this solution satisfies $\sum_j \Vert (X_j^R)_t \Vert_\infty\leq R$, 
it will also be a solution to the original equation 
\eqref{sde} involving $DV$.
Furthermore, if we start with some initial 
data $X_0$ with $\sum \Vert (X_j)_0\Vert_\infty <R$, solutions to
\eqref{sde} always exist for small time (at least up until the time
that the operator norm of the solution exceeds $R$).
Thus we may consider a solution up to the
time that its norm reaches some fixed constant $M$ (with the intent 
to show that this
time
is infinite).
By free It\^o calculus
\begin{eqnarray*}
dX_t.X_t&=&
2X_t .dS_t-DV(X_t).X_t dt+2dt\\
&=&2X_t .dS_t -(DV(X_t)-DV(0)).X_t dt-DV(0).X_t dt
+2dt
\end{eqnarray*}
Therefore, for all $s\ge 0$,
\begin{eqnarray*}
e^{c s} X_s.X_s&=&
X_0.X_0 -2\int_0^s e^{ct}DV(0).X_t dt
-\int_0^s e^{ct}[(DV(X_t)-DV(0)).X_t -c X_t.X_t]dt\\
&&+2\int_0^s e^{ct}dt+2\int_0^s e^{ct}.dS_t \\
&\le& X_0.X_0 -\int_0^s e^{ct}DV(0).X_t dt+ 2c^{-1}e^{
cs}+2\int_0^s e^{ct}.dS_t
\end{eqnarray*}
where the last inequality holds in the
sense of operator and we relied on our hypothesis of $(c,M)$ convexity.
 Since also $e^{c s} X_s.X_s$
is a non negative operator,
we deduce that 
\begin{eqnarray*}
A_s&:=&\|X_s.X_s\|_\infty\\
&\le& e^{-cs}A_0+2c^{-1}+
2\left\|\int_0^s e^{c(t-s)}
X_t .dS_t\right\|_\infty+
\|DV(0).DV(0)\|_\infty^{\frac 12}
\int_0^s e^{c(t-s)} A_t^{\frac 12} dt
 \end{eqnarray*}
By Theorem 3.2.1 of \cite{BS},
we know that the free analog of the Burkh\"older-Davis inequality
for integrals with respect to free
Brownian motion holds for the $L^p$ norm even with $p=\infty$.
More precisely, the following estimate holds:
$$\left\|\int_0^s e^{c(t-s)}
X_t .dS_t\right\|_\infty\le 2\sqrt{2} 
\left(\int_0^s e^{2c(t-s)} A_t dt\right)^{\frac{1}{2}}
\le  2\sqrt{2} 
\left(\int_0^s e^{c(t-s)} A_t dt\right)^{\frac{1}{2}}
.$$
Moreover, by the Cauchy-Schwarz inequality, we obtain the bound 
$$\left(\int_0^s e^{c(t-s)} A_t^{\frac 12} dt\right)^2
= \left( \int_0^s e^{\half c(t-s)} \cdot e^{
\half c(t-s)} A^{\frac{1}{2}}_t dt\right)^2\le
\frac{1}{c} \int_0^s e^{c(t-s)} A_t dt.$$
Hence,  we get the inequality (with
$C=\frac{4}{c}(\|DV(0).DV(0)\|_\infty^{\frac 12}+4\sqrt{2})^2$
and  since $A_0\le b^2 m$)
\begin{eqnarray*}
A_s^2&\le& 4e^{-cs}m^2b^4 +2^4c^{-2}+
C\int_0^s e^{c(t-s)} A_t dt\\
&\le 
& 4e^{-cs}m^2b^4+2^4c^{-2}+
\frac{C}{2}\int_0^s e^{c(t-s)}(B^{-1} A_t^2 +B) dt\\
&\le 
& 4e^{-cs}m^2b^4+2^4c^{-2}+CB(2c)^{-1}+
CB^{-1}\int_0^s e^{c(t-s)} A_t^2 dt
\end{eqnarray*}
where $B$ is any positive
constant (we just used that for all $x$, $2x\le B+ B^{-1} x^2$). 

We now use Gronwall's lemma (or simply iterate in
the above inequality)
to deduce, with $C'=2^4c^{-2}+CB(c)^{-1}$, that 
$$A_s^2 e^{cs}
\le [4m^2b^4+C'e^{cs}] +\frac{C}{B}\int_0^s 
e^{\frac{C}{B}(s-u)} [4m^2b^4+C'e^{cu}] du$$
We now take  $B$ large enough so that $c> \frac{C}{B}$
to conclude that
$$A_s^2 
\le 8m^2b^4e^{(\frac{C}{B}-c)s} +C'\frac{c}{c-\frac{C}{B}}$$
If we now choose $B_0^2= (2^4c^{-2}+CB(c)^{-1})\frac{c}{c-\frac{C}{B}}$
for $B=2C/c$, we have shown  that:
\begin{itemize} 
\item $A_s^2$ stays bounded 
by $8m^2b^4 +B_0^2$. So if $8m^2b^4<M^2- B_0^2$,
$A_s$ always stays bounded by $M$.
\item As $s$ goes to infinity, $\limsup_{s\ra\infty} A_s$
is bounded by $B_0$.
\item We can choose $b>B_0$ as long as $M> B_0\sqrt{1+8m^2 B^2_0}:=M_0$. 
\item $X_t\in C^*(Z,S_q:q\geq 0)$.
\end{itemize}
This concludes the proof.
\hfill$\square$

We now consider the solutions of  Lemma \ref{exist} starting with different
initial data.

\begin{theo}\label{existtheo}
Let $M_0$, $B_0$ and $b$ be as in Lemma~\ref{exist}, and assume that
$M\geq M_0$, and that $Z$ is an 
$m$-tuple of operators with $\Vert Z\Vert_\infty < b.$ 
Consider the unique solutions $X_t^Z$, $X_t^0$ to the free SDE
\begin{equation}\label{freeSDE}
dX_t=dS_t-\half DV(X_t)dt
\end{equation}
with initial conditions $X_0^Z=Z$, $X_0^0 = 0$.  
Then 
\begin{enumerate}
\item $\Vert X_t^Z - X_t^0 \Vert_\infty \to 0\qquad\textrm{as $t\to\infty$}.$

\item The law of $((X_t^Z)^*,X_t^Z)$
converges to a stationary 
law  $\mu_V\in\cxmnsa'$
which satisfies for all $P\in\cxmnsa$ and all $i\in\{1,\ldots,m\}$

\begin{equation}\label{SDtot1}
\sum_{i=1}^m \mu_V\otimes\mu_V\left(  (\partial_i+\partial_{i,*})
(D_i+D_{i,*}) P\right)=\sum_{i=1}^m\mu_V\left( D_i V D_iP+ (D_i V)^*
D_{i,*} P\right).
\end{equation}
Moreover, for any $k\in\N$,
any $i\in \{1,\ldots, m\}$,
\begin{equation}\label{compact}
\mu_V\left((X_iX_i^*)^k\right)\le B_0^{2k}.
\end{equation}
Any law $\nu\in\cxmnsa'$ of variables bounded
in operator norm by $b$
which satisfies \eqref{SDtot1} 
equals $\mu$ on $\cxmnsa$. 
 
\item
The restriction $\tau_V= \mu_V|_\cxm$ 
of $\mu_V$ to $\cxm$ satisfies for all polynomials $P\in\cxm$ 
\begin{equation}\label{SD}
\sum_{i=1}^m\tau_V\otimes\tau_V(  \partial_i D_i P)=\sum_{i=1}^m \tau_V(D_i V D_i P).
\end{equation}

Any law $\nu\in\cxmnsa'$ of variables bounded in operator norm by $b$
whose retriction $\nu|_\cxm$
satisfies \eqref{SD} is such that $\nu|_\cxm= \tau_V$. 
\end{enumerate}
\end{theo}

\begin{rmk}
(a) Note that if $Z$ is an $m$-tuple of variables, norm bounded by $b$
having the stationary law $\mu$, then $X^Z_t$ 
is a stationary process (with law given at all times by
$\mu$) and 
 $X_t^Z-X_t^0$ converges to zero in operator norm
as $t$ goes to infinity by the first point of the theorem. 

(b) Recall that a law $\mu$ of $m$ non-commutative
variables has conjugate variables $(\xi_i)_{1\le i\le m}$
if and only if for any $P\in\cxmnsa$,
$$\mu\otimes\mu(\partial_i P)=\mu(\xi_i P).$$
Taking the adjoint, we find that
we must also have 
$$\mu\otimes\mu( \partial_{i,*}P)=\mu(\xi_i^* P).$$
Taking $P=D_iQ$ in the first equality, $P=D_{i,*} Q$ in
the second and summing the resulting equalities yields
$$\mu\otimes\mu((\partial_i D_i + \partial_{i,*}  D_{i,*})Q)
= \mu( \xi_i D_i Q+ \xi_i^* D_{i,*} Q)$$
which differs from 
\eqref{SDtot1} by the terms
$\mu\otimes\mu(\partial_i D_{i,*} Q+ \partial_{i,*}  D_{i}Q )$. Hence, 
\eqref{SDtot1} is not equivalent with the fact that $\mu$
 has conjuguate the variable $(D_iV)_{1\le i\le m}$ in the case that 
$V$ is not self-adjoint.   In the self-adjoint case,  because of our
convention, $Q$ depends only on $X_1,\ldots,X_n$ and so
the terms involving $D_{i,*}$ and $\partial_{i,*}$ 
are equal to zero. In that case, \eqref{SDtot1} is compatible with the 
condition that the conjugate variables are equal to the cyclic gradient
of $V$.
\end{rmk}

{\bf Proof.}
Consider two solutions $X_t^Y$, $X_t^Z$ with initial data $Y$ and $Z$,
respectively, and assume that $\Vert Y\Vert_\infty, \Vert Z\Vert_\infty 
\leq m$.
Then
$$d(X^Z_t-X^Y_t)=-\half [DV(X_t^Z)-DV(X^Y_t)]dt$$
Since by Lemma~\ref{exist}, the operator norms of $X^Z_t$ and $X^Y_t$
 stay bounded by $M$ for all $t$ and $V$ is 
$(c,M)$ convex, we find that
\begin{eqnarray*}
d(X^Z_t-X^Y_t).(X^Z_t-X^Y_t)  &\le& -[DV(X_t^Z)-DV(X^Y_t)].(X^Z_t-X^Y_t)
dt\\
&\le&-c (X^Z_t-X^Y_t).(X^Z_t-X^Y_t)dt
\end{eqnarray*}
where the inequality again holds in the
sense of operators.
This implies that
\begin{equation}\label{control1}
\|(X^Z_t-X^Y_t).(X^Z_t-X^Y_t)\|_\infty\le e^{-ct} \|Z-Y\|^2_\infty 
\end{equation}
and so 
$$\lim_{t\ra\infty} \|(X^Z_t-X^Y_t).(X^Z_t-X^Y_t)\|_\infty=0.$$
In particular, we can take $Y=0$
and we have proved that all diffusion solutions starting
from different initial data norm-bounded by $b$ will 
asymptotically be the same as $X^0_t$.

As $t$ gets large, $X^0_t$ is bounded by
$B_0<b$ according to Lemma \ref{exist},
and so we can choose $Y= \tilde X_t^0$
(with $\tilde X$ constructed as $X$ but with some
 free Brownian motion $\tilde S$ on $[0,s]$
and the increments of $S$ on $[s,s+t]$ )
  to
deduce 
\begin{equation}\label{cauchy}
\|(X^0_t-\tilde X^0_{t+s}). (X^0_t-\tilde X^0_{t+s}) \|_\infty \le e^{-ct} M.
\end{equation}
Since $\tilde X^0_{t+s}$
has the same law as $X^0_{t+s}$
we conclude that the joint law of $X^0_t,(X^0_t)^*$ 
converges as $t$ goes to infinity.
We denote this limit by $\mu$.
 Clearly, $\mu$ is a stationary
law for the diffusion (since $X^0_{\infty+s}$ 
and $X^0_{\infty+t}$ have the same law).
We now consider $X^\infty_t$ to
be the solution of the free SDE starting from $X$ 
so that $(X, X^*)$ has law  $\mu$.
Noting that we have 
$$dX^\infty_t=dS_t-\half DV(X_t^\infty) dt \qquad
d(X^\infty_t)^*=dS_t-\half  (DV(X_t^\infty))^* dt$$
and applying 
 free It\^o's calculus \cite{BS},
we have that for any  $P\in\cxmnsa$,
\begin{multline*}
0=\partial_t \phi(P(X^\infty_t, (X^\infty_t)^* ))\\
=\phi\otimes\phi(\frac{1}{2} \sum_{i=1}^m (\partial_{i,*}+
\partial_i)(D_{i,*}+ D_i) P(X^\infty_t, (X^\infty_t)^*  )) \\-\half \phi(
DV(X^\infty_t).DP(X^\infty_t)+(DV(X^\infty_t))^*. D_*P(X^\infty_t ,(X^\infty_t)^* )
) )\end{multline*}
 so that for all $*$-polynomials $P\in\cxmnsa$,
$\mu_V$ must satisfy \eqref{SDtot1}. $\tau_V$, the restriction of $\mu_V$ to
$\cxm$ must then satisfy \eqref{SD}.

The uniqueness of solutions
to this equation  is simply due 
to the fact that if we run 
the process from $Z$ bounded uniformly by $b$ with law  $\nu$
satisfying \eqref{SDtot1}, the law of $X^Z_t$ must
be stationary,
but also converging to $\mu_V$ by the previous argument.
Hence it must be equal to $\mu_V$.  The same argument applies in the case that 
we only care about the restriction to $\cxm$
of some law $\nu$ of variables $Z$ bounded by $b$,
 since if it satisfies \eqref{SD}, the process $X^Z_t$
will be such that $\partial_t \phi( P(X^Z_t))=0$
for all $P\in \cxm$ and so the law of $X^Z_t$ restricted
to $\cxm$ will be stationary (here 
we use the fact that $DV(X)$ depends only on $X$ and not on
 its adjoint when speaking of the evolution of the
law of $X^Z_t$ restricted to $\cxm$). Since it also converges 
to $\tau_V$, we obtain the desired equality.
\hfill$\square$

We shall see in section \ref{analytic}, 
Corollary \ref{corSD},
that in fact $\mu$ not only satisfies
\eqref{SD} but actually for all $i\in\{1,\ldots,m\}$,
all $P\in\cxm$,
$$\mu\otimes \mu (\partial_i P)=\mu(D_i V P).$$
In other words, at least when $\mu$
is the law of self-adjoint operators,  the  conjugate variables 
of $\mu$ are in
the cyclic gradient space.

\section{Analyticity  of the solution 
to Schwinger-Dyson equation and discussion
around phase transition}\label{analytic}
We show in this section
that on the domain where $V$ stays
$(c,M)$-convex for some
$c>0$ and $M\ge M_0$
as in Lemma \ref{exist},
the law $\tau_V$ will depend
analytically on the parameters 
of $V$.

\begin{lem}\label{analy}
Let $V=V_{\beta}=\sum_{i=1}^n \beta_i q_i$ be a polynomial, where
${\bf \beta}=(\beta_i)_{1\le i\le n}$ are (complex) parameters 
and $(q_i)_{1\le i\le n}$ are monomials. For $c,M$ positive real
numbers and $M\ge M_0$ with $M_0$ as in Lemma \ref{exist}
and Theorem
\ref{existtheo},
let $T(c,M
)\subset \C^n$ 
be the interior 
of the subset of parameters ${\bf \beta}=(\beta_i)_{1\le i\le n}$
for which $V$ is $(c,M)$-convex.
Let $\mu_\beta=\mu_{V_\beta}$ be the 
unique stationary
measure of Theorem \ref{existtheo}
and $\tau_\beta$ the law of $(X_1,\cdots,X_m)$ under
$\mu_\beta$.  

Then  
for any polynomial $P\in \cxm$, the map ${\bf \beta}\in T(c,M) \ra 
\tau_\beta(P)$ is analytic. 
\end{lem}
Note that $T(c,M)$ is non empty 
as soon as the set of monomials 
$(q_i, 1\le i\le n)$ 
contains $(X_i^2, 1\le i\le m)$.
Indeed, if we set $V(X)=\sum_{i=1}^m \beta_i X_i^2 +\sum_{i=m+1}^n \beta_i
q_i(X)$, we always have that for $X,Y$ uniformly bounded
by $M$, 
$$\left(\sum_{i=m+1}^n \beta_i
D q_i(X)- \sum_{i=m+1}^n \beta_i
D q_i(Y)\right).(X-Y)\le C(M)\max_{m+1\le i\le n}|\beta_i| (X-Y).(X-Y)$$
with a universal constant $C(M)$ which only depends on $M$ and the $(q_i)_{1\le i\le n}$.
Hence, if $\beta_i>0$ for $i\in\{1,\cdots,m\}$,
$$(DV(X)-DV(Y)).(X-Y)\ge [\min_{1\le i\le m}\beta_i-C(M)\max_{m+1\le i\le n}
|\beta_i|
]  (X-Y).(X-Y).$$
Therefore any set of parameters $(\beta_i)_{1\le i\le n}$
such that $$\min_{1\le i\le m}\beta_i-C(M)\max_{m+1\le i\le n}
|\beta_i|\ge c$$
will be such that $V$ is $(c,M)$-convex.

\noindent
{\bf Proof.}

We denote $(X^{\bf \beta}_t)_{t\ge 0}$ the solution 
of \eqref{freeSDE} with potential $V=V_{\bf \beta}$
and starting from the null operator.
We shall show that $\beta \ra X^{\beta}_t$ expands as 
a  sum of uniformly bounded operators.
More precisely, we fix $\beta$ in the interior
of $ T(c,M)$
 and find 
a family $X^{(k_1,\ldots,k_n)}_t, k_i\in\N, 1\le i\le n$
of  operator-valued  processes such
that for  $\eta\in \C^n$,
$|\beta-\eta|:=\max_{1\le i\le n}|\beta_i-\eta_i|$
small enough,
\begin{equation}\label{res1}
X^{\bf \eta}_t=X^{\bf \beta}_t+
\sum_{\stackrel{k_1,\ldots,k_n\in\N^n}{\sum k_i\ge 1}}
 \prod_{i=1}^n  (\eta_i-\beta_i)^{k_i} X^{(
k_1,\ldots,k_n
)}_t
\end{equation}
Moreover,  $X^{
(
k_1,\ldots,k_n
)
}_t$
are operator-valued  processes such that there exists 
a constant $C$ which only depends on $c,M$ and the degree 
of $V$ so that 
\begin{equation}\label{res12}\sup_{t\in\R^+}
\| X^{
(
k_1,\ldots,k_n
)
}_t\|_\infty \le C^{\sum k_i }.\end{equation}
Finally the distribution of $\left(X^{
(
k_1,\ldots,k_n
)
}_t\right)_{k_1,\ldots, k_n\in\N^n}$ 
converges (in the sense of finite marginals, i.e., on polynomials involving
only a finite number of the $\left(X^{(k_1,\ldots,k_n)}_t\right)_{k_1,\ldots,k_n\in \N^n}$)
towards the law of  $\left(X^{(k_1,\ldots,k_n)}_\infty\right)_ {k_1,\ldots, k_n\in\N^n}$
as $t$ goes to infinity. 

Let us conclude the proof of the lemma
assuming \eqref{res1} and \eqref{res12}. \eqref{res1} and \eqref{res12}
entail that for all polynomial functions 
$P\in\cxm$, for all $t\ge 0$,
$$\beta \ra \phi(P( X^{\bf \beta}_t))$$
is analytic in the interior of $T(c,M)$
since it implies that for  $\beta\in
T(c,M)$, 
$$\phi_\beta(P(X^{\bf \eta}_t))=\phi     \left
( P(
X^{\bf \beta}_t+\sum_{k_1,\ldots,k_n}
 \prod_{i=1}^n  (\eta_i-\beta_i)^{k_i} X^{(k_1,\ldots,k_n)}_t)
\right)$$
for all $\eta$ in the domain 
$B(C,\beta)=\{|\eta-\beta|< 1/C\}$
which does not depend on the time parameter $t\ge 0$.
Note also that $\tau(P( X^{\bf \eta}_t))$ is uniformly bounded
independently of $t\in \R^+$ since $C$ does not depend
on $t$. We know that $(X^{\bf \eta}_t)$
 converges as $t$ goes to infinity towards 
$$\tau_{\beta}(P):=\phi( P(
X^{\bf \beta}_\infty+\sum_{k_1,\ldots,k_n}
 \prod_{i=1}^n  (\eta_i-\beta_i)^{k_i} X^{(k_1,\ldots,k_n)}_\infty))$$
(note here that convergence of 
the $ X^{(k_1,\ldots,k_n)}$ in the sense of finite marginals 
is sufficient since
 $\sum_{k_1,\ldots,k_n:\sum k_i\ge
K}
 \prod_{i=1}^n  (\eta_i-\beta_i)^{k_i} X^{k_1,\ldots,k_n}_\infty$
goes uniformly to zero as $K$ goes to infinity).
 But then the limit has to depend analytically on $\beta$
(as a limit of uniformly bounded  functions which are analytic
on
a fixed domain).
This proves the claim.

We now prove \eqref{res1} and \eqref{res12}.
We first check that $\beta\ra X^{\beta}_t$ 
is of  class $\Ca^\infty$, then that 
it 
is in fact an
entire function and bound uniformly its radius
of convergence. Finally, we prove that 
$t\ra ( X^{(k_1,\ldots,k_n)}_t)_{k_1,\cdots, k_n\in\N^n}$ 
converges.

{\it Step 1: $\beta\in T(c,M)\ra X^{\beta}_t$ 
is of class $\Ca^\infty$ for all $t\ge 0$.}

Let us  study  the first order
differentiability, and 
 first check that ${\bf \beta}\ra X^{{\bf \beta}}_t$
is continuous.
In fact, if $1_p(i)=0$ for $i\neq p$
and $1_p(p)=1$, we write
\begin{eqnarray*}
X^{{\bf \beta}+\e 1_p}_t -X^{{\bf \beta}}_t
&=& - \half
\int_0^t [ DV_{\beta+\e 1_p}(X^{{\bf \beta}+\e 1_p}_s) 
- DV_{\beta}(X^{{\bf \beta}}_s)]ds\\
&=&-
\half\int_0^t 
\int_0^ 1 \partial DV_{ \beta+\e 1_p} ((1-\alpha) X^{{\bf \beta}}_s
+\alpha  X^{{\bf \beta}+\e 1_p}_s)\sharp [X^{{\bf \beta}+\e 1_p}_s -X^{{\bf \beta}}_s]
d\alpha ds\\
&& -\half\int_0^t  [DV_{\bf \beta+\e 1_p }-DV_{\bf \beta}]( 
X^{{\bf \beta}}_t).
dt
\end{eqnarray*}
We find that if ${\bf \beta}+\e 1_p$
and $\beta$ both belong to $T(c,M)$ 
so that $X^{\bf \beta}$ and $X^{\beta+\e 1_p}$
stay  uniformly bounded by $M$, that
$$\| X^{{\bf \beta}+\e 1_p}_t -X^{{\bf \beta}}_t\|_\infty
\le C \int_0^t \| X^{{\bf \beta}+\e 1_p}_s -X^{{\bf \beta}}_s\|_\infty
+C t \e$$
where $C$ only depends on $(c,M)$. 
So Gronwall's lemma shows that for any time $t\ge 0$, 
there exists a finite $C(t)$
($C(t)$ is uniformly bounded on compacts)
 so that 
\begin{equation}\label{bound11}\| X^{{\bf \beta}+\e 1_p}_t -X^{{\bf \beta}}_t\|_\infty
\le C(t)\e.\end{equation}
This suggests that ${\bf \beta}\ra X^{{\bf \beta}}_t$
is in fact differentiable. To prove this  point, let
us introduce the candidate for the corresponding gradient;
we 
 define $(\nabla_{\bf \beta} X^{\beta,i}_t)_{1\le i\le m}$
to be the $m$-tuple 
of operators valued processes solution of
$$d \nabla_{\bf \beta}^p X^{ \beta,i}_t= -\sum_{j=1}^m
\partial_j D_i V(X^{\bf \beta}_t) \sharp \nabla_{\bf \beta}^p X^{ \beta,j}_t
dt + (d_\beta^p D_i V_{\bf \beta})(X^{\bf \beta}_t) dt
\qquad 
\nabla_{\bf \beta}^p X^{\beta,i}_0\equiv 0
$$
where $d_\beta$ is the standard gradient with
respect to the parameters  $\beta$
(so $d_\beta^j D_i V_{\bf \beta}=D_i q_j$ for $j\in\{1,\ldots, n\}$).
Here $p$ is any integer in $\{1,\ldots,n\}$.
There is a unique solution to this
equation (since
it is a linear differential equation as 
 $X^{\bf \beta}$
is given).
By the same type
of argument as above, we now prove that 
for all $t\ge 0$ there exists
$C(t)$ finite so that
\begin{equation}\label{bound12}
\| X^{{\bf \beta}+\e 1_p}_t -X^{{\bf \beta}}_t-\e 
\nabla_{\bf \beta}^p X^{\bf \eta,i}_t\|_\infty
\le C(t)\e^2
\end{equation}
Indeed, 
if we let 
$$Y_t^{\e,i}:=X^{{\bf \beta}+\e 1_p, i}_t -X^{{\bf \beta},i}_t-\e 
\nabla_{\bf \beta}^p X^{ \beta, i}_t$$
we find that 

\begin{multline*}
Y_t^{\e}
=-
\int_0^t 
\int_0^ 1 \partial DV_{\bf \beta+\e 1_p} ((1-\alpha) X^{{\bf \beta}}_s
+\alpha  X^{{\bf \beta}+\e 1_p}_s)\sharp [Y_s^{\e}] ds\\
 - \int_0^t 
\int_0^ 1[ \partial  DV_{\bf \beta+\e 1_p}
 ((1-\alpha) X^{{\bf \beta}}_s
+\alpha  X^{{\bf \beta}+\e 1_p}_s)-\partial DV_{\bf \beta+\e 1_p}
 ( X^{{\bf \beta}}_s)] \sharp [
X^{{\bf \beta}+\e 1_p}_s -X^{{\bf \beta}}_s] ds\\
  -\int_0^t  [DV_{\bf \beta+\e 1_p }-DV_{\bf \beta}-\e d_\beta^p D_i V_{\bf \beta} ]( 
X^{{\bf \beta}}_t)
dt
\end{multline*}
Using \eqref{bound11} 
we find that there exists a finite constant $C=C(M)$ 
such that
\begin{eqnarray*}
\max_{1
\le i\le m}\|Y_t^{\e,i}\|&\le&
 C\int_0^t \max_{1
\le i\le m}
\|Y_s^{\e,i}
\| ds + C\e^2 \int_0^t C(s)^2 ds
\end{eqnarray*}
and so Gronwall's lemma gives \eqref{bound12}
(note here that $\max_{1
\le i\le m}
\|Y_t^{\e,i}
\|$ is finite for all $\e>0$ 
so that $\beta+\e1_p\in T(c,M)$). 

This shows that ${\bf \beta}\ra X^{{\bf \beta}}_t$
is differentiable for all $t$
with first derivative $ \nabla_{\bf \beta} X^{ \beta}_t$.
We can continue in the same spirit 
to show that $ \nabla_{\bf \beta}^p X^{\beta,i}_t$
is differentiable and by induction, 
we find that ${\bf \beta}\ra X^{{\bf \beta}}_t$
is of class $\Ca^\infty$ in the interior of  $T(c,M)$.
We next bound uniformly all its derivatives.

{\it Step 2: $\beta\in T(c,M)\ra X^\beta_t$
is analytic  for all $t\ge 0$ }.

To this end, we can write $X^\beta_t$
as a formal series in $\eta$
in a small ball around $\beta\in T(c,M)$
$$X^\eta_t= \sum_{(k)=(k_1,\ldots, k_n)\in\N^n}  
\prod_{1\le p\le n}(\eta_p-\beta_p)^{k_p} X^{(k)}_t$$
with $X^{
(
0,\ldots,0
)
}_t=X^\beta_t$. 
Indeed, the coefficients of this series are obtained
by differentiating the $\Ca^\infty$ operator $X^\eta_t$
and  $ X^{(k)}_t=(k_1!\cdots k_n!)^{-1} \partial_{\eta_1}^{k_1}\cdots 
\partial_{\eta_n}^{k_n} X^\eta_t
|_{\eta=\beta}$.
We write $D_iV_\beta(X)=\sum_{j=1}^n \beta_j D_i q_j
=\sum_{j=1}^{Dn}\tilde \beta_j  q_{ij}$
where $D_iq_j=\sum_{l=1}^{D} q_{i,l+(j-1)D}$
is the decomposition of $D_iq_j$ as a sum of 
at most $D$ monomials. We denote  $q_{ij}=\prod^{\ra}_{1\le p\le d_{ij}}
X_{l_{ij}^p}$
with $l_{ij}^p\in \{1,\ldots,m\}$.
 Moreover, we have   $\tilde\eta_j=\eta_{[j/D]}$.
Plugging this formal series into 
$$dX^\eta_t= dS_t -\half DV_\eta(X^\eta_t)dt$$
we find that $X^{
(
k_1,\ldots, k_n
)
}_.$
satisfy, for $\sum k_i\ge 1$,  the following induction
relation
\begin{eqnarray}
dX^{(k),j}_t&=&-\half\sum_{i=1}^{nD}
\sum_{\sum_{p=1}^{d_{ji}}
 k_r^p= k_r-1_{r=[i/D]}}
\prod^{\ra}_{1\le p\le d_{ji
}}
X^{(k^p),l_{ji}^p}
 dt\nonumber\\
&&-\half\sum_{i=1}^{nD}\tilde \beta_i
\sum_{\stackrel{\sum_{p=1}^{d_{ij}
}
k_r^p= k_r-1_{r=[i/D]}}{ ( k^p)\neq (k)\forall
p}}
\prod^{\ra}_{1\le p\le d_{ji}}
X^{(k^p),l_{ji}^p} dt\nonumber\\
&&-\half\sum_{l=1}^m \partial_l D_j V_\beta(X^\beta_t)\sharp X^{(k), l} dt
\label{lkj}
\end{eqnarray}
Above, the sum over indices $k$ such that
$k_r-1_{r=[i/D]}=-1$ is simply empty. 
Using the convexity of $V_\beta$, we now get a uniform 
bound by considering $X^{(k)}_t.X^{(k)}_t=
\sum_{j=1}^m X^{(k),j}_t.X^{(k),j}_t$;
\begin{gather}
dX^{(k)}_t.X^{(k)}_t\le -\sum_{j=1}^m\sum_{i=1}^{nD}
\sum_{\stackrel{\sum_{p=1}^{d_{ji}} k_r^p=} {k_r-1_{r=[i/D]}}}
\left[ \big(\prod^{\ra}_{1\le p\le d_{ji}}
X^{(k^p),l_{ji}^p}\big)^*X^{(k),j}_t+
(X^{(k),j}_t)^*\prod^{\ra}_{1\le p\le d_{ji}}
X^{(k^p),l_{ji}^p}\right]
 dt\nonumber\\
-\sum_{j=1}^m\sum_{i=1}^{nD}\tilde \eta_i
\sum_{p=1}^{d_{ji}}\sum_{ \stackrel{k_r^p= k_r}{( k^p)\neq (k)\forall
p}}\left[ \big(\prod^{\ra}_{1\le p\le d_{ji}}
X^{(k^p),l_{ji}^p}\big)^*X^{(k),j}_t+
(X^{(k),j}_t)^*\prod^{\ra}_{1\le p\le d_{ji}}
X^{(k^p),l_{ji}^p}\right] dt\nonumber\\
- cX^{(k)}_t.X^{(k)}_t dt \label{lkj22}\\
\nonumber
\end{gather}
where we have used that
by convexity,
for any operators $X,Z$ bounded by $M$,
$$\sum_{i=1}^m \sum_{j=1}^m ( \partial_j D_i V(X)\sharp Z_j).Z_i
\ge c Z.Z.
$$
Note that in the right hand side 
of \eqref{lkj22}, all the $X^{(k^p),l_{ji}^p}$
are such that $\sum k^p_i<\sum k_i$.
Hence, we can deduce by
induction
that $A^{(k)}_t:=
\max_{1\le j\le m}\|X^{(k),j}_t\|_\infty$ is bounded 
for all $(k)$ and uniformly on
compact sets of the time variable  $t$.  
Indeed, we proved it in Theorem \ref{existtheo}
for $(k)=(0,\ldots,0)$.
Let us assume it is true for $(k)$ with $\sum k_i\le K-1$
and let us prove it remains true; we simply use
$$\left\|\prod^{\ra}_{1\le p\le d_{ji}}
X^{(k^p),l_{ji}^p}.X^{(k),j}_t\right\|_\infty
\le \prod_{1\le p\le d_{ij}} A^{(k_p)}_t A^{(k)}_t
\le B \prod_{1\le p\le d_{ij}} (A^{(k_p)}_t)^2 +B^{-1} ( A^{(k)}_t)^2.$$
Choosing $B$ such that
$$ 2B^{-1}
\sum_{j=1}^m \sum_{i=1}^{nD}(1+|\beta_i|) d_{ji}\sum_{\sum k_r^p\le k_r}1
-c<0$$
allows to bound $A^{(k)}_t$ uniformly 
on compact sets by our induction hypothesis.
We next show that this bound
can be taken uniformly
on the time variable.
To this end we first consider $Y^{(k)}_t=\sqrt{X^{(k)}_t.X^{(k)}_t}$
and deduce from
\eqref{lkj22} that
\begin{eqnarray}
dY^{(k)}_t&\le &m\sum_{i=1}^{nD}\ \ 
\sum_{\sum_{p=1}^{D} {k_r^p=}{ k_r-1_{r=[i/D]}}} \ \ 
\prod_{1\le p\le D}
A^{(k^p)}_t dt\nonumber\\
&&+mD \sum_{i=1}^{n}|\beta_i|
\sum_{ \stackrel{\sum_{p=1}^{D} k_r^p= k_r}{ ( k^p)\neq (k)\forall
p}}\prod_{1\le p\le D} A^{(k^p)}_t dt -c Y^{(k)}_t
dt
\label{lkj3}
\end{eqnarray}
where we eventually added terms (by taking $d_{ij}=D$)
 which can be done
if we assume that $A^{(0)}\ge 1$
which we can always do.
Therefore, for 
$\sum k_i\ge 1$ 
, since $X^{(k),j}_0=0$,
we obtain the bound
\begin{eqnarray}
A^{(k)}_t
&\le& mD \sum_{i=1}^n \int_0^t e^{-c(t-s)}
\sum_{\sum_{p=1}^{D} k_r^p= k_r-1_{r=i}}
\prod_{1\le p\le D}A^{(k_p)}_s
 ds\nonumber\\
&& +mD \sum_{i=1}^n|\beta_i| \int_0^t e^{-c(t-s)}
\sum_{\stackrel{\sum_{p=1}^{D} k_r^p= k_r}{ ( k_p)\neq (k)\forall
p}}
\prod_{1\le p\le D}A^{(k_p)}_s
 ds.\label{chouette}
\end{eqnarray}
Note that the right hand side 
of \eqref{chouette} depends 
only on $A^{(l)}$
for $\sum l_i\le \sum k_i-1$
(since $ k_r^p\le k_r$  for all $r$
but $(k^p)\neq (k)$).
Since $A^{(0)}$ is uniformly bounded (as we proved
that $X^\beta_t $ is uniformly bounded 
for $\beta\in T(c,M)$  and $c>0$), we
deduce by induction that
$A^{(k)}:=\sup_{t\ge 0}A^{(k)}_t$
is finite  and satisfy
the induction bound, with $|\beta|=\sum_{i=1}^n |\beta_i|$,
\begin{eqnarray}
A^{(k)}&\le& mD c^{-1}\sum_{i=1}^n 
\sum_{\sum_{p=1}^{D} k_r^p= k_r-1_{r=i}}
\prod_{1\le p\le D}A^{(k_p)}\nonumber\\
&& +mDc^{-1} |\beta| 
\sum_{\stackrel{\sum_{p=1}^{D} k_r^p= k_r}{ ( k_p)\neq (k)\forall
p}}
\prod_{1\le p\le D}A^{(k_p)}\label{chouette2}
\end{eqnarray}
We rewrite this inequality,
since $A^{(k)}$ is obviously finite,
as
\begin{eqnarray}
A^{(k)}&\le& \frac{ mD c^{-1}}{1+
mD^2 c^{-1}|\beta| (A^{(0)})^{D-1}}\Big\{
\sum_{i=1}^n 
\sum_{\sum_{p=1}^{D} k_r^p= k_r-1_{r=i}}
\prod_{1\le p\le D}A^{(k_p)}\nonumber\\
&& +mDc^{-1}|\beta|
\sum_{\sum_{p=1}^{D} k_r^p= k_r}
\prod_{1\le p\le D}A^{(k_p)}\Big\}\label{chouette3}
\end{eqnarray}
where we added the term $mDc^{-1}|\beta| D (A^{(0)})^{D-1} A^{(k)}$
to the last sum. 
We now want to show that there is
a finite $C=C(\beta)$ so that
$$A^{(k)}\le C^{\sum k_i}.$$
To this end we borrow the idea
of majorizing sequences as developed by Cartan \cite{cartan},
chapter VII.
It goes as follows here.
We consider the polynomial in
one variable given by 
$$\tilde V_{\kappa}(x)= C \kappa x^D\quad\mbox{ with }
C:= \frac{ mD c^{-1}}{1+
mD^2 c^{-1}|\beta| (A^{(0)})^{D-1}}.$$
and the 
the equation
$$ x-x_0= \tilde V_{\kappa}(x) -\tilde V_{
|\beta|}(x_0)$$
with $x_0=A^{(0)}$.
We claim that for $\kappa$ in a neighborhood
of $
|\beta|$,
the solution $x_\kappa$ in the neighborhood of $x_0$
of this equation is analytic in $\kappa$.
Indeed, by the implicit function theorem,
we only need to check that $$1\neq \partial_x \tilde
 V_{|\beta|}(x_0)=\frac{ mD c^{-1}|\beta|}{1+
mD^2 c^{-1}|\beta| (A^{(0)})^{D-1}}  D (A^{(0)})^{D-1}$$
which is always true since $A^{(0)}$ is finite. 
This implies that there exists a finite constant $C=C( |\beta|)$
such that
$$| \partial_{\kappa}^k
x_{\kappa}|_{\kappa=|\beta|}|
\le k! C^{k}.$$
But now, $x^{(k)}:=(k!)^{-1}\partial_{\kappa}^{k}
x_{\kappa}|_{\kappa=|\beta|}$
satisfies the induction relation
$$
x^{(k)}=C \sum_{ \sum_{i=1}^D  k_i=k-1} \prod_{1\le i\le D} x^{(k_i)}
+C|\beta| 
\sum_{\sum_{i=1}^D  k_i=k} \prod_{1\le i\le D} x^{(k_i)}$$
which implies
$$(1-C  |\beta|
Dx_0^{D-1})x^{(k)}= C \sum_{ \sum_{i=1}^D  k_i=k-1} \prod_{1\le i\le D} x^{(k_i)}
+C|\beta| 
\sum_{\sum_{i=1}^D  k_i=k, k_i<k} \prod_{1\le i\le D} x^{(k_i)}.$$
Since $1-C |\beta| Dx_0^{D-1}= 1/(1+mDc^{-1}|\beta|
 D x_0^{D-1})>0$, we conclude by induction that
 $x^{(k)} \ge 0$ for all $k$  (note that $x^{(0)}=A^{(0)}>0$).
But then, comparing \eqref{chouette2} and
the above inequality, we also prove by induction that 
$$A^{(k)}\le x^{(\sum k_i)}\le C^{\sum k_i}$$
which therefore
gives  the desired bound for the $A^{(k)}$.
Hence we have proved \eqref{res12}.

{Step 3: \it Convergence in law of of $t\ra (X^{(k)}_t)_{(k)}$
for all $(k)$}.

As we have noticed
above, the equations for $X^{(k)}$ are of the form
$$dX^{(k)}_t= -\partial DV_\beta( X^\beta_t)\sharp X^{(k)}_t
dt +P_{(k)} (X^{(l)}_t, \sum l_i\le \sum k_i -1) dt$$
with some polynomial functions $P_{(k)}$. 
Therefore, if we denote $X^{(k),Z}_t$
the solution of this equation starting from
$X^{(k)}_0=Z^{(k)}$,
we get from the convexity of $V$ and by induction 
over $\sum k_i$
that
$$\sup_{\sum k_i\le K}\|X^{(k),Z}_t -X^{(k),0}_t\|_\infty
\le C(K) e^{-t}$$
with some finite constant $C(K)$.
Hence, we can start from $Z^{(k)}=X^{(k)}_{s}$
to see that 
since $(X^{(k),Z}_t)_{k\in\N^n}$ has
the same law as $(X^{(k),0}_{t+s})_{k\in\N^n}$,
the law of $(X^{(k),0}_t,,  \sum k_i\le K)$
converges. 
\hfill$\square$

We now relate the previous result 
with the absence of phase transition
for the generating function of colored planar maps.
In \cite{GM}, the following  strong version
of Schwinger-Dyson equation was considered; it requires  that for
 all   polynomials $P\in\cxm$,  
\begin{equation}\label{stronSDd}
\tau(D_i V P)=\tau\otimes\tau(\partial_i P),\, 1\le i\le m.
\end{equation}
It was shown that if $V(X_1,\ldots,X_m)=W_{(q_i,\beta_i)_{1\le i\le
n}}=\frac{1}{2}\sum_{i=1}^m X_i^2
+\sum_{i=m+1}^{n} \beta_i q_i$, 
 there exists
 a unique solution  $\Ma\in \cxm'$   under the condition that
$$|\tau(X_{i_1}\cdots X_{i_k})|\le R^k$$
for all $k$ and some finite $R$, provided 
 the  $\beta_i$'s  are small enough. We denote 
$\tau_{W_{(q_i,\beta_i)_{1\le i\le
n}}}
$ this solution.    Note that this
solution was  not a priori  the law of non-commuting 
variables, except in the case where $V$ is self-adjoint, but just
an element of $\cxm'$. In particular, E. Maurel Segala and one of
the authors always restricted 
to polynomials in the letters $(X_1,\ldots,X_m)$ and did not
consider their adjoints.

Moreover, for all monomial $P\in\cxm$,
$$\tau_{W_{(q_i,\beta_i)_{1\le i\le
n}}}(P)=
\Ma_{(q_i,\beta_i)_{1\le i\le
n}} (P)=\sum_{k_1,\ldots,k_n\in\N^n}\prod_{1\le i\le n}
\frac{(-\beta_i)^{k_i}}{k_i!}M_0( (q_i,k_i), 1\le i\le n,(P,1))$$
with $M_0( (q_i,k_i),(P,1))$ the number of planar maps
with $k_i$ stars of type $q_i$ for
$1\le i\le n$ and one star of type $P$ (we refer the reader to \cite{GM},
section 2, 
for a complete description of the numbers $M_0( (q_i,k_i),1\le i\le n,(P,1))$).

We now claim
\begin{theo}\label{corSD}
(a) The generating function 
$$(\beta_i)_{1\le i\le n} \in\C^n\ra \Ma_{(q_i,\beta_i)_{1\le i\le
n}} (P):=\sum_{k_1,\ldots,k_n}\prod_{1\le i\le n}
\frac{(-\beta_i)^{k_i}}{k_i!}M( (q_i,k_i),1\le i\le n,(P,1)),$$
which is an absolutely convergent series
for $\sum_{i=1}^n |\beta_i|$
small enough, 
extends analytically in the interior of the domain 
where $W_{(q_i,\beta_i)_{1\le i\le
n} }
(X_1,\ldots,X_m)=\half \sum_{i=1}^m X_i^2
+\sum_{i=1}^n \beta_i q_i$ is $(c,M)$-convex
for some $c>0$ and
$M\ge M_0(c, \|DV(0).DV(0)\|_\infty)$
the constant of Theorem \ref{existtheo}. The extension of $
\Ma_{(q_i,\beta_i)_{1\le i\le
n} }(P)$
is equal to $\tau_{W_{(q_i,\beta_i)_{1\le i\le
n} }}$, the restriction to $\cxm$
of the invariant measure of Theorem
\ref{existtheo}. \\
(b) Assume that $V$ is $(c,M)$-convex
with $M\ge M_0$ of Theorem \ref{existtheo}. The invariant distribution  $\mu_{V}$ of Theorem
\ref{existtheo} not only satisfies \eqref{SD}
but its strong version in the sense that $\tau_V=\mu_V|_\cxm$
is such that
\begin{equation}\label{stronSDdd}
\tau_V(D_i V P)=\tau_V\otimes\tau_V(\partial_i P),\, 1\le i\le m
\end{equation}
for all polynomials $P\in \cxm$.
\end{theo}
Hence, the first point of the above theorem shows that 
the breaking of
analyticity (or phase transition)
of the map enumeration
can not take place when $W_{(q_i,\beta_i)_{1\le i\le
n}}$ is $(c,M)$-convex.

\noindent
{\bf Proof.} By \cite{GM},
if we consider the case where $W_{(q_i,\beta_i)_{1\le i\le
n}}$
is self-adjoint, we know that $\Ma_{(q_i,\beta_i)_{1\le i\le
n}}$
is the law of self-adjoint
operators which
are uniformly bounded by $R$ ($R=R(\beta)$ going
to  $2$
as $\beta$ goes to zero)  when the
$\beta$'s are small enough. As a consequence,
$\Ma_{(q_i,\beta_i)_{1\le i\le
n}}$ must coincide with $\tau_{1,(\beta_i)_{1\le i\le n}
}=\mu_{1,(\beta_i)_{1\le i\le n}
}
|_\cxm$ where we put $\mu_{1,(\beta_i)_{1\le i\le n}
}:=\mu_{W_{(q_i,\beta_i)_{1\le i\le
n}
}}$ 
since $\Ma_{(q_i,\beta_i)_{1\le i\le
n}}$ satisfies \eqref{SD} with potential $W_{(q_i,\beta_i)_{1\le i\le
n}}$ (by Theorem \ref{existtheo}).
In particular, $\tau_{1,(\beta_i)_{1\le i\le n}}$
must satisfy \eqref{stronSDdd}. 

We now show that we can remove the assumption 
that $W_{(q_i,\beta_i)_{1\le i\le
n}}$  is self-adjoint. We denote by
$*$ the involution 
$(z X_{i_1}\cdots X_{i_k})^*= \bar
z X_{i_k}\cdots X_{i_1}$
so that $W_{(q_i,\beta_i)_{1\le i\le
n}}$ is self adjoint iff $W_{(q_i,\beta_i)_{1\le i\le
n}}
=W_{(q_i,\beta_i)_{1\le i\le
n}}^*$.
We can always write $W_{(q_i,\beta_i)_{1\le i\le
n}}$ in the form 
$W_{(q_i,\beta_i)_{1\le i\le
n}, (q_i^*,\beta_i')_{1\le i\le
n}
} = \half \sum X_i^2 
+\half  \sum \beta_i q_i +\half \sum \beta_i' q_i^*$. We denote
by $\mu_{1,\beta,\beta'}$
the invariant measure of
Theorem \ref{existtheo} corresponding
to such a potential.
The situation $V=V^*$
corresponds to $\beta_i'=\bar \beta_i$.
In that case we have shown that
\begin{equation}\label{qwe}
\tau_{1,(\beta_i)_{1\le i\le n},(\beta_i')_{1\le i\le n}}(D_iV P)
-\tau_{1,(\beta_i)_{1\le i\le n},(\beta_i')_{1\le i\le n}}\otimes
\tau_{1,(\beta_i)_{1\le i\le n},(\beta_i')_{1\le i\le n}}(\partial_i P)=0.
\end{equation}
But recall that if an analytic function of two variables $x,y$ 
is null on $\Lambda=\{x=\bar y, |x|\le \e\}$
for some $\e>0$, then this function must
vanish on its full domain of analyticity since $\Lambda$ is totally real.
Hence, inside the domain of analyticity of $\beta,\beta'\ra
\tau_{1,(\beta_i)_{1\le i\le n},(\beta_i')_{1\le i\le n}}$ (i.e the domain where
$(\beta_i)_{1\le i\le n},(\beta_i')_{1\le i\le n}
\ra
\tau_{1,(\beta_i)_{1\le i\le n},(\beta_i')_{1\le i\le n}}
(P)$ is analytic for all $P\in\cxm$),
  \eqref{qwe} is always true.
We can now  remove the artificial parameters $(\beta_i')_{1\le i\le n}$ 
to claim that $\tau_{1,(\beta_i)_{1\le i\le n}}$ always  satisfy \eqref{stronSDdd}
on its domain of analyticity.
We finally can remove the condition that $\beta_i=1$ 
for $1\le i\le m$ as follows. We take $V$ $(c,M)$-convex with $M\ge M_0$.
 Note first that 
 if we consider $V_\alpha=\frac{\alpha}{2} \sum_{i=1}^m
X_i^2 + V$ for some $\alpha>0$, then the result
still holds since by uniqueness of
the solution to Schwinger-Dyson equation, we have the scaling property
$$\mu_{V_\alpha}=d_\alpha\sharp \mu_{
\frac{1}{2} \sum X_i^2 + d_\alpha\sharp V
}$$
with $d_\alpha$ the dilatation 
$d_\alpha \sharp \mu(P)= 
\mu(P(\frac{X_1}{\sqrt{\alpha}}, \cdots, \frac{X_m}{\sqrt{\alpha}}))$. 
As $d_\alpha\sharp V$ is always $(0,M)$ convex,
 $\frac{1}{2} \sum X_i^2 + d_\alpha\sharp V$ satisfies the above hypotheses and so 
 $\tau_{\frac{\alpha}{2} \sum_{i=1}^m
X_i^2 + \sum_{i=1}^m \beta_i q_i}$ always satisfies \eqref{stronSDdd}. 
Finally, we can let $\alpha$ going to zero since $V_\alpha$
is $(c,M)$-convex for all $\alpha\ge 0$
and so $\alpha\ra \mu_{V_\alpha}(P)$ is analytic
and
thus continuous when $\alpha$ goes to zero.

As a conclusion, we have seen that $\tau_{V_\beta}$
satisfies \eqref{qwe}  on  the domain of analyticity
which contain by Lemma \ref{analy}
all the sets $T(c,M)$, $c>0, M>M_0(c)$.

On the other hand,
we also have that $\tau_{1,\beta}$
must agree with $\Ma_{\beta}$, the generating function
of maps, for all $\beta$ small enough
since $\tau_{1,\beta}$ satisfies
 \eqref{stronSDdd} and is the law of bounded operators
(and it was proved that there is at most one such solution,
 the generating function of maps, in \cite{GM}).
We hence conclude that $\beta\ra \Ma_{\beta}(P)$
extends analytically to the domain
of analyticity of $\tau_{1,\beta}$,
which contains all $\beta$ such that $(1,\beta)$
belongs to $T(c,M)$
for $c>0$ and $M\ge M_0(c)$.

\hfill$\square$

\section{Connectivity of the support and properties of associated $C^*$ and 
von Neumann algebras}\label{connect}

Throughout the rest of the paper, we shall assume that $V$ is a 
$(c,M)$ convex potential with $M>M_0$ (so that the hypothesis of
Theorem~\ref{existtheo} holds).

In this section, 
we denote by $\mu$ the unique stationary law for the free stochastic 
differential equation (\ref{freeSDE}) with drift 
$DV$, satisfying \eqref{SD}, where $S_q$ is a free Brownian motion.   Lastly,
$Z$ will denote some fixed $m$-tuple of operators, free from $S_q:q\geq 0$,
having law $\mu$, and satisfying $\Vert Z\Vert_\infty < b$.

The main results of this section concern properties of the $C^*$-algebra
generated by the $m$-tuple $Z$ with the prescibed law $\mu_V$.  We show
that this $C^*$-algebra is exact \cite{wassermann}, projectionless and that the associated
von Neumann algebra has the Haagerup approximation property 
\cite{haagerup:compact} and admits
and embedding into the ultrapower of the hyperfinite II$_1$ factor.  These
properties are shared by (and in fact, in most cases, derived from those of) 
the $C^*$-algebra generated by $m$ semicircular systems.

One of the most interesting open
problems  in operator algebras is a 
question due to Connes of whether any 
tracial state has finite approximation
in the sense that there exists a norm-bounded sequence
of $N\times N$ matrices $(A_1^N,\ldots, A_m^N)$ 
such that for all *-polynomial function $P$,
$$\lim_{N\ra \infty}
\frac{1}{N}\tr(P( A_1^N,\ldots, A_m^N))=\tau(P).$$
When $m=1$, this question is settled by Birkhoff's
theorem, but the question is still
open when $m\ge 2$.  We prove that the laws $\mu_V$ have finite-dimensional
approximations.

\subsection{Approximation of $Z$ by elements from $C^*(S_q: q\geq 0)$.}
We first show that, with respect to operator norm, $C^*(S_q : q\geq 0)$ 
$\varepsilon$-contains a variable with law $\mu$.
\begin{cor} \label{corNormApprox}
Within the hypothesis and notations of Theorem
\ref{existtheo},  
for any $\varepsilon > 0$, there exists a Brownian motion $S_s : s\geq 0$ free from $Z$,
and elements $X'\in C^*(S_s : s\geq 0)$, $X\in C^*(Z,S_s:s\geq 0)$  so that 
$X$ has the given stationary law $\mu$ and 
$\Vert X - X'\Vert_\infty \leq \varepsilon$.
\end{cor}
{\bf Proof.}
Let $X_t$ and $X_t^Z$ be two solutions to \eqref{sde} with initial data
$X_0=0$, $X_0^Z=Z$.
By Theorem \ref{existtheo},
$X_t\in C^*(S_s : s\leq t)$ approximates in operator norm the
stationary
process $X_t^Z$ 
with marginal distribution $\mu$.  Since
by Lemma~\ref{exist}, 
$X_t^Z\in C^*(Z,S_q:q\geq 0)$, we may take $X'=X_t$, $X=X_t^Z$ for 
large enough $t$.  
\hfill$\square$

\begin{theo}
Let $Z$ be any $m$-tuple of $b$-bounded variables
with the unique law $\mu$ satisfying
\eqref{SD}.
(a) The algebra $C^*(Z)$ has no non-trivial projections.  
(b) The spectrum of any non-commutative
*-polynomial $P$ in the $m$-tuple $Z$ is connected (in the case that 
$P(Z)$ is normal, this means that the support of its spectral measure is 
connected).  
(c) If $P$ is any polynomial in $Z$ whose value is self-adjoint, then 
the probability measure given by the law of $P(Z)$ has connected support.
\end{theo}
{\bf Proof.} The $C^*$-algebra $\mathcal{A}=C^*(S_s, s\geq 0)$  can be identified
with the $C^*$-algebra generated by semicircular operators 
$s(f): f\in L^2(\mathbb R; \mathbb R)$ where $f\mapsto s(f)$ denotes
the free Gaussian functor \cite{BS}.  It is well  known
that this $C^*$-algebra is isomorphic to the 
infinite reduced free product
$$\mathcal{A}\cong (C[-1,1],\mu) * \cdots *(C[-1,1],\mu)$$ where $\mu$ denotes
the semicircular measure.  The algebra $(C[-1,1],\mu)$ can be 
unitally embedded in a trace-preserving way into the 
group $C^*$-algebra $C^*(\mathbb Z)\cong C(\mathbb T)$ taken with 
its canonical group trace $\tau$.   Indeed, if
$u\in C^*(\mathbb Z)$ denotes the group generator, i.e. $u=\exp(2\pi i\theta)
\in C(\mathbb T)$, then $u+u^*$ generates a copy of $C[-1,1]$,
and the restriction of $\tau$ to $C^*(u+u^*)$ is the arcsine law.
Hence for a suitable continuous function $f$, $f(u+u^*)$ has as its 
distribution the semicircle law, and we can embed 
$(C[-1,1],\mu)$ by sending its generator, multiplication by $x$,
to $f(u+u^*)$. It follows
that $$\mathcal{A}\subset C^*(\mathbb Z) *_\textrm{red} 
\cdots *_\textrm{red} C^*(\mathbb Z) 
\cong C_\textrm{red}^*({\mathbb F}_\infty).$$
By the results of \cite{PV} (see also \cite{HST} for a random 
matrix proof), $C^*({\mathbb F}_\infty)$ has no 
non-trivial projections.  Thus $\mathcal{A}$ has
no non-trivial projections.

Suppose now that $Y\in \mathcal{A}$.  Then the spectrum 
$\sigma(Y)$ must be connected.  We sketch the argument, which 
can be found in standard $C^*$-algebra literature (see e.g.
\cite[Proposition 4.6.2 on p. 28]{blackadar:book}).  If $\sigma(Y)=K_1 \cup K_2$
with $K_1\cap K_2 =\emptyset$ and both $K_1, K_2$ non-empty, for any 
contour $\gamma \subset \mathbb C$ that contains $K_1$ but not $K_2$ and
does not intersect $\sigma(Y)$, the integral
$$ E = \frac{1}{2\pi i} \int_\gamma (z-Y)^{-1}  dz$$ belongs to
$\mathcal{A}$,
being a norm limit of Riemann sums.
Moreover, $E^2=E$, $E(1-E)=0$ but 
 $E\neq 0$, $E\neq 1$.  
   But 
 $E+E^* - 1$ is invertible since $(E+E^*-1)^2 = 1 + (E-E^*)(E-E^*)^*
\geq 1$, and thus $$P = E(E+E^*-1)^{-1}$$ is a nonzero self-adjoint
projection (since $(1-E)E=E^*(1-E^*)=0$). 
Clearly, 
$P\neq 0$ and $(1-E)P = 0$ so $P\neq 1$.  
Hence $\mathcal{A}$ 
would have a non-trivial
projection, a contradiction.

Since $X_t^0\in \mathcal{A}$, it follows that the spectrum of 
any  $*$-polynomial $Y$ in this $m$-tuple is connected.
Since $Z$ and $X^Z_t$ have the same laws at all $t$, $C^*(Z)$ and 
$C^*(X^Z_t)$ are isomorphic.  So 
if the support of the spectrum of $P(Z)$ is disconnected,
then the support of the spectrum of  $P(X^Z_t)$ is disconnected, for all $t$.
Because  $X_t^0 $ converges to $ X_t^Z$
in operator norm as $t\to\infty$, we find that also spectrum of 
$P(X_t^0)$ must be disconnected for large enough $t$.  Since the 
trace-state on $C^*(S_s:s\geq 0)$  is faithful, it follows that for some
$t$ large enough, $C^*(S_s: s\geq 0)$ contains an element 
with disconnected spectrum.  But as we saw before, this is impossible.
Thus we have proved (b).

Because  any non-trivial projection has disconnected spectrum, (a) follows.

Lastly, as  the GNS construction for the trace-state on $C^*(Z)$ is
faithful (by construction of $C^*(Z)$) and the state is tracial, it follows
that the trace vector in the representation is cyclic for both 
$C^*(Z)$ and its commutant.  But this implies that the vector is also 
separating, so that the trace-state on $C^*(Z)$ is faithful.  Thus the support
of the law of any self-adjoint operator $X\in C^*(Z)$ is exactly its spectrum.
This implies (c).
\hfill$\square$

\subsection{Exactness and the Haagerup property.}

We recall that a $C^*$-algebra $A$ is called exact (cf.
\cite{wassermann} and references therein) if there exists
a faithful $*$-representation $\pi : A\to B(H)$ with the following property.
For any finite subset $F\subset A$ and any $\varepsilon > 0$ there exists
a finite-dimensional matrix algebra $D$ and 
unital completely positive maps $\theta : A\to D$, $\eta: D\to B(H)$ so that
$$\Vert \eta(\theta(x)) - \pi (x) \Vert_\infty < \varepsilon,\qquad
\forall x\in F.$$
It turns out that this property is equivalent to the statement that the
functor $\otimes_\textrm{min} A$ of taking the minimal tensor product with $A$ 
is exact (i.e., takes exact sequences to exact sequences).  Exactness
is an important approximation property for a $C^*$-algebra.

Another important approximation property, this time for a von Neumann
algebra, is the Haagerup property.
A von Neumann algebra $M$ with a trace $\tau$ is said to have the Haagerup
property \cite{haagerup:compact} if there exists a sequence of 
completely-positive maps $\Phi_n : M\to M$, which are unital and trace-preserving, so that the associated maps $\Phi_n : L^2(M,\tau)\to L^2(M,\tau)$ are 
compact and converge to $1$ strongly: $$
\Vert \Phi_n (x) - x \Vert_2 \to 0,\qquad n\to \infty$$ for all 
$x\in M$.  For a discrete group von Neumann algebra $L(\Gamma)$, the 
Haagerup property is 
equivalent to Gromov's a-T-menability of the group $\Gamma$ (i.e., to the 
existence of a cocycle $c: \Gamma \to H$ with values in some unitary 
representation $H$ of $\Gamma$, so that the map $\gamma \mapsto \Vert c(\gamma)\Vert_H$ is proper).  As was shown by Haagerup, free groups have the Haagerup
property.
\begin{theo}
Let $Z$ be any $b$-bounded $m$-tuple with 
the unique law $\mu$ satisfying \eqref{SD}.
Then $C^*(Z)$ is exact and $W^*(Z)$ has the Haagerup approximation property.
\end{theo}
{\bf Proof.}
As in Corollary~\ref{corNormApprox}, 
let $X^Z_t$ be the solution to the free SDE (\ref{freeSDE}) starting with $Z$
and $X^0_t$ be the solution starting with zero.  Thus $X^Z_t \in C^*(Z, S_q:
q\in [0,+\infty))$ and $X_t \in C^*(S_q : q\in [0,+\infty))$, where $S_q$ is a 
free Brownian motion, free from $C^*(Z)$.  Moreover, we have by Theorem~\ref{existtheo}
that $\Vert X_t^Z - X_t^0\Vert_\infty \to 0$ as $t\to \infty$.

Let $E: C^*(Z, S_q : q\in [0,+\infty))\to C^*(S_q : q\in [0,+\infty))$ be
the conditional expectation coming from the fact that $Z$ and $S_q, q\in [0,+\infty)$  are 
freely independent.  
Then 
\begin{equation}
\label{eq:distBrM}
\Vert E(Z_t - X_t )\Vert_\infty \leq \Vert Z_t -X_t \Vert_\infty\to 0,
\end{equation}
so that, since $E(X_t)=X_t$, $$\Vert E(Z_t)-X_t \Vert_\infty \to 0.$$
Since $\mu$ is stationary, $C^*(X^Z_t)\cong C^*(Z)$; let $\pi_t : C^*(Z)\to 
C^*(X^Z_t)$ be this isomorphism.  Let $H=L^2(W^*(Z,S_q:q\in [0,+\infty)))$. Then 
as a module over $W^*(Z)$, $H$ is infinite-dimensional, and as a module over
$W^*(X^Z_t)$, it is at most infinite-dimensional.  Hence there exists an injective unital $*$-homomorphism
$\Theta_t : B(H) \to B(H)$ with the property that $\Theta( \pi_t (x)) = x$ for all 
$t\geq 0$ and $x\in W^*(Z)$ (and thus for $x\in C^*(Z)$).  

Let now $F\subset C^*(Z)$ be a finite subset and $\varepsilon > 0$.  Then one can find
$*$-polynomials $P_x : x\in F$ so that $$\Vert P_x(Z)- x \Vert_\infty <\varepsilon /9 .$$  Hence
for $t$ sufficiently large, we may assume that $$\Vert E(P_x(X_t^Z)) - P_x(X_t^Z)\Vert_\infty
<\varepsilon / 9.$$  Since $\pi_t(P_x(Z))= P_x(X_t^Z)$, we conclude that 
$$\Vert E(\pi_t (x)) - \pi_t(x) \Vert_\infty <\varepsilon / 3,\qquad \forall x\in F.$$
Now for all $x\in F$, $E(\pi_t(x))\in C^*(S_q : q \in [0,+\infty))$, and hence we can find
a finite-dimensional $C^*$-algebra $A$ and unital completely positive maps $\eta : C^*(S_q : 
q\in [0,+\infty))\to A$, $\psi : A \to B(H)$, so that $$\Vert \psi \circ \eta (y) - y \Vert_\infty < \varepsilon / 3$$ for all $y\in \pi_t(F)$. 
 Consider now the unital completely positive maps
$$\alpha =  \eta \circ E\circ \pi_t : C^*(Z)\to A$$ and 
$$\beta = \Theta \circ \psi : A\to B(H).$$
Then $$\Vert \beta\circ\alpha (x) - x \Vert_\infty =  
\Vert \eta(\pi_t(x)) - \pi_t (x)\Vert_\infty < 
\varepsilon$$ for all $x\in F$.  Thus $C^*(Z)$ is exact.

We now turn to the Haagerup property, where we adapt a proof from
\cite{Br}. 
Consider the map $\Phi_t : W^*(X^Z_t) \to W^* (X^Z_t)$ which is obtained as the composition
$$\Phi_t = E_{W^*(X^Z_t)} \circ \Psi_t \circ E_{W^*(S_q:q\in[0,+\infty))}, $$
where $\Psi_t$ are unital trace-preserving completely positive maps on 
$W^*(S_q : q\in [0,+\infty))$ so that $\Psi_t$ are compact on $L^2$ and
$\Vert \Psi_t (x) -x \Vert_2 \to 0$ for all $x\in W^*(S_q : q\in [0,+\infty))$.
Then $\Phi_t$ are unital 
trace-preserving completely-positive maps on $W^*(X_t^Z)\cong W^*(Z)$,
and because of (\ref{eq:distBrM}), one has that $\Vert \Phi_t(x)-x \Vert_2 \to 0$ for all
$x\in W^*(Z)$.  On the other hand, each $\Phi_t$, viewed as a map on 
$L^2 (W^*(Z, S_q : q\in [0,+\infty)))$, is compact (since $\Psi_t$ is compact),
and therefore the restriction of $\Phi_t$ to $L^2(W^*(X_t^Z))$ is also compact.
Thus  $W^*(Z)$ has the Haagerup property.
\hfill$\square$
\subsection{Finite dimensional approximation.}\label{finiteap}
\subsubsection{$R^\omega$ embeddability for self-adjoint potentials.}

In this section, we improve on the results of \cite{GM}
by showing that if we
set
$$\bar\mu^N_V(P)=\int \hat\mu_N(P) d\mu^N_V(A_1,\ldots,A_m))$$
with $$\hat\mu_N(P)=
\frac{1}{N}\tr( P(A_1,\ldots,A_m))$$
and 
$$d\mu^N_V(A_1,\ldots,A_m))=\frac{1}{Z_N^V} 1_{\|A_i\|_\infty\le M}
  e^{-N\tr(V(A_1,\ldots,A_m))}
dA_1\cdots dA_m, $$
then $\bar\mu^N_V$ converges towards $\tau_V$ 
for any self-adjoint locally strictly convex potential $V$.
In \cite{GM}, a similar result was proved
when $V=cX.X+W$ with
$W$ a 'small' enough polynomial.

\begin{theo}\label{conjuguate}
For all $c>0$,
there exists $B_0<\infty$ and $M_0<\infty$
so that
for any self-adjoint polynomial $V$
which is $(c,M)$-convex
with $M>M_0$
there exists a unique  law  $\tau$   of $m$ self-adjoint variables 
such that for all $i\in\{1,\ldots,m\}$,
all polynomial $P$, 
\begin{equation}\label{SDstrong}
\tau\otimes \tau (\partial_i P)=\tau(D_i V P)
\end{equation}
and  such  that $\tau(X_i^{2k})\le B_0^{2k}$.
Moreover, $\bar\mu^N_V$ converges towards $\tau$ 
and therefore 
$\tau$ has finite approximation.  In particular, if $Z$ has law 
$\tau$, then $W^*(Z)$ can be embedded into the ultrapower of the hyperfinite
II$_1$ factor.
\end{theo}
{\bf Proof.}
We can follow the lines of \cite{GM} Theorem 3.5
to see that $d\mu^N_V(A_1,\ldots,A_m))$ has a density with
respect to the Gaussian law $\prod e^{-\frac{Nc}{2} X_i^2} dX_i/Z_N^c$
which is log-concave. This insures that we can use the Brascamp
Lieb inequality which in turn allows us to show
that the random matrices under the above Gibbs measures
stay bounded in norm by some $B_0\ll M$ with
overwhelming probability. As a consequence, we can perform
an infinitesimal change of variables  $X_i\ra X_i+N^{-1} P_i$
with $P_i$ a self-adjoint polynomial in $\|A\|_\infty\le B_0$
and the  null function
outside   a ball of radius strictly smaller than $M$.
This shows that any almost sure limit points of
$\hat\mu_N$ under $\mu^N_V$ are laws of variables bounded 
by $B_0$ and satisfy \eqref{SDstrong} with the potential $V$  (see \cite{GM} for details).
Now, if they satisfy  \eqref{SDstrong}
they also satisfy  \eqref{SD} (take $P=D_iP$
and sum the equalities) and so
by Theorem \ref{existtheo}, 
there exists at most one such solution $\mu_{V}=\tau_{V}$.
Thus, $\hat\mu_N$ converges almost surely and
therefore in expectation
towards $\tau_{V}$. Hence $\tau_{V}$ has finite approximation.
\hfill$\square$

\subsubsection{$R^\omega$ embeddability for non-self adjoint potentials.}
One can give a proof of embeddability of $W^*(Z)$ into the ultrapower of the hyperfinite
II$_1$ factor, based directly on Corollary~\ref{corNormApprox}.  This proof works for arbitrary
$(c,M)$-convex polynomials (without the self-adjoint assumption).  Indeed, 
because of this corollary and with its notations,  if  $\varepsilon >  0$ is given, and
$F$ is a finite collection of $*$-polynomials, then there exists an $X'\in C^*(S_q:q\geq 0)$
with the property that 
$$|\tau(P(Z)) - \tau(P(X'))|<\varepsilon,\qquad \forall P\in F$$
(here $\tau$ denotes the free product trace-state on $C^*(Z,S_q:q\geq 0)$).  
Now, 
$C^*(S_q: q\geq 0)$ is generated by an infinite free semicircular family $\hat{S}_j : j=1,2,\ldots$.
One can clearly assume that in fact $X'\in C^*(\hat{S}_1,\ldots,\hat{S}_K)$ for some large enough $K$
(since one can replace $X'$ with $E_K(X')$ for $K$ large enough, 
where $E_K : C^*(\hat{S}_i : i=1,2\ldots)
\to C^*(\hat{S}_1,\ldots,\hat{S}_K)$ is the canonical conditional expectation).  Thus, by approximating
$X'$ with a polynomial in $\hat{S}_1,\ldots,\hat{S}_K$, we may assume that $X'_k$ is a polynomial
$Q_k$
in $\hat{S}_1,\ldots,\hat{S}_K$.
Since $\hat{S}_1,\ldots,\hat{S}_K$ are free semicircular variables, their law has finite
approximations.  Thus 
for $N$ sufficiently large, one can find a $K$-tuple of $N\times N$ self-adjoint matrices
$A_1,\ldots,A_K$ whose law approximates that of $\hat{S}_1,\ldots,\hat{S}_K$ so well that
the $m$-tuple of matrices 
$B=(B_1,\ldots,B_m) =(Q_1(A_1,\ldots,A_K),\ldots,Q_m(A_1,\ldots,A_K))$ would have the 
property that  $$|\tau (P(Z))-\frac{1}{N}\tr
 (P(B))|<\varepsilon,\qquad \forall P\in F.$$

\section{Free
entropy}
For the remainder of the paper we shall assume that $V$ is 
$(c,M)$-convex and self-adjoint.

In this section we
show that for tracial state with conjugate variables
given as the cyclic gradient of a self-adjoint  $(c,M)$-convex potential
the microstate entropy
is the same whether it
is defined by a limsup or a liminf. 
\begin{theo}\label{entropy} Let
$c>0$ and $V$ be a self-adjoint $(c,M)$-convex
potential with $M>M_0$. 
Let $\tau=\tau_V$ be as in Theorem  \ref{conjuguate}.
Let $\Gamma(\tau,\e,R,k)$ be the microstates
\begin{multline*}
\Gamma(\tau,\e,R,k)=\{ X_1,\cdots,X_m:\|X_i\|_\infty\le R, 
| \frac{1}{N} \tr(P(X_1,\cdots,X_m))-\tau(P)|\le \e,
\\ \mbox{ for all monomial of degree less than }k\}
\end{multline*}
and let $\mbox{vol}$ denote
the volume on the space of $m$ $N\ts N$ Hermitian
matrices.

Then 
\begin{eqnarray*}
\chi(\tau)&=& \limsup_{\e\downarrow 0,k, R\ra \infty}
\limsup_{N\ra\infty} \frac{1}{N^2}\log
{\mbox vol}(\Gamma(\tau,\e,R,k)) + \frac{m}{2}\log N\\
&=&\liminf_{\e\downarrow 0,k, R\ra \infty}
\liminf_{N\ra\infty} \frac{1}{N^2}\log
{\mbox vol}(\Gamma(\tau,\e,R,k)) +\frac{m}{2}\log N
\end{eqnarray*}
and $\chi(\tau)>-\infty$.
\end{theo}
{\bf Proof.}
Note that
\begin{eqnarray}
{\mbox vol}(\Gamma(\tau,\e,R,k))&=&\int 1_{A_1,\cdots, A_m\in \Gamma(\tau,\e,R,k)} dA_1\cdots dA_m\\
&\approx&
e^{N^2\tau(V)}
\int_{\Gamma(\tau,\e,R,k)} e^{-N^2\mun(V)}dA_1\cdots dA_m\nonumber\\
&\approx&e^{N^2\tau(V)}
\mu^N_V(\Gamma(\tau,\e,R,k))\int e^{-N^2\mun(V)}dA_1\cdots dA_m
\label{poi}\\
\nonumber
\end{eqnarray}
with $\mu_N^V$ the Gibbs measure considered in
the proof of Lemma \ref{conjuguate}. Here, we used
the notation $A(N,\e)\approx B(N,\e)$
when
$$\lim_{\e\downarrow 0}\lim_{N\uparrow 0}
\frac{1}{N^2}\log \frac{A(N,\e)}{B(N,\e)}=1.$$
Since $\mu^N_V(\Gamma(\tau,\e,R,k))\to 1$ 
by the proof of Theorem \ref{existtheo}, we only need to estimate
the quantity $\int e^{-N^2\mun(V)}dA_1\cdots dA_m$.
To do that we write $V=W+\frac{c}{2}X.X$
with $W$ a $(0,M)$-convex potential
so that
\begin{eqnarray*}
\partial_t\frac{1}{N^2}\log
 \int e^{-N^2\mun(tW+\frac{c}{2}X.X)}dA_1\cdots dA_m
&=&-\bar \mu^N_{tW+\frac{c}{2}X.X}(W)
\end{eqnarray*}
By the proof of Theorem  \ref{conjuguate},
we also see that for all $t\in [0,1]$
$$\lim_{N\ra\infty}
\bar \mu^N_{tW+\frac{c}{2}X.X}(W)=\tau_{tW+\frac{c}{2}X.X}(W)$$
(note that $tW+\frac{c}{2}X.X$ stays self-adjoint $(c,M)$-convex).
Since everything
stays bounded,  and since the limit 
$$\lim_{N\ra\infty} \frac{1}{N^2} \log 
\int e^{-N^2 \mun(\frac{c}{2}X.X)} dA_1\cdots dA_m +\frac{m}{2}\log N
$$ is a finite constant $F(c)$, 
we conclude by bounded convergence theorem
that
$$\lim_{N\ra\infty}
\frac{1}{N^2}\log
 \int e^{-N^2\mun(W+\frac{c}{2}X.X)}dA_1\cdots dA_m + \frac{m}{2}\log N
=F(c)-\int_0^1 \tau_{tW+\frac{c}{2}X.X}(W)dt
$$
which shows by \eqref{poi} that 
$\chi$ can be defined either
by limsup or liminf
and both are equal to
$$\chi(\tau)=F(c)-\int_0^1 \tau_{tW+\frac{c}{2}X.X}(W)dt
+\tau(V)$$

\begin{cor}
Let $\tau$ be as in Theorem  \ref{conjuguate}.
Then the von Neumann algebra generated by an $m$-tuple $Z$ with law $\tau$ 
is a factor, does not have property $\Gamma$, is prime and has no Cartan subalgerbas.
\end{cor}

Here we of course use the fact that $\chi(\tau)>-\infty$ and 
that for an $m$-tuple $Z$, 
$\chi(Z)>-\infty$ implies that $W^*(Z)$ is a factor, non-$\Gamma$, has no Cartan subalgebras
\cite{dvv:entropy3} and is prime \cite{ge:prime}.

\section{Norm convergence}
In section \ref{finiteap}
it was shown that solutions to the Schwinger-Dyson
equations \eqref{SDstrong}
are weak limits of finite dimensional approximations.
Namely if we let
$$\bar\mu^N_V(P)=\int \hat\mu_N(P) d\mu^N_V(A_1,\ldots,A_m)$$
with $$\hat\mu_N(P)=
\frac{1}{N}\tr( P(A_1,\ldots,A_m))$$
and 
$$d\mu^N_V(A_1,\ldots,A_m))=\frac{1}{Z_N^V} 1_{\|A_i\|\le M}
  e^{-N\tr(V(A_1,\ldots,A_m))}
dA_1\cdots dA_m,$$
we saw that $\bar\mu^N(P)$ converges to the unique solution $\tau$ 
of 
$$\tau_V(D_i V P)=\tau_V\otimes\tau_V\left( \partial_i P\right)$$
for all $i\in\{1,\ldots, m\}$ and $P\in \cxm$.
We here show that this convergence holds 
in norm but for simplicity restrict ourselves
to potentials which are uniformly convex (and not only locally convex).
\begin{lem}
Let $(X_1,\ldots,X_m)$ be a m-tuple
of non-commutative variables with law $\tau_V$
and let $(X_1^N,\ldots,X_m^N)$ be random matrices with
law $\mu^N_V$. Assume that $V$ is self-adjoint 
 $(c,\infty)$-convex
for some $c>0$. 
Then, for any polynomial $P\in\cxm$,
$$\lim_{N\ra\infty}
\| P(X_1^N,\ldots,X_m^N)\|_\infty
=\| P(X_1,\ldots,X_m)\|_\infty\,\mbox{a.s.}$$
\end{lem}
This result generalizes the work of Haagerup and
S. Thorjornsen \cite{HT}
where it was proved for $V(X_1,\ldots, X_m)=\frac{1}{2}
\sum_{i=1}^m X_i^2$. Our result actually relies on theirs.

\noindent {\bf Proof.}
The idea is to use the approximation by processes,
and hence  the fact that processes are well approximated by polynomials 
of independent Wigner matrices and then use 
\cite{HT} to conclude that the norm of the latter converge
to the limit.

{\it Step 1: Matrix valued diffusions and
convergence to the stationary
process.}

Consider the diffusion
with values in the set of Hermitian matrices
\begin{equation}\label{eqdiff}
dX^{N,Z}_t= dH^N_t-\half DV(X^{N,Z}_t) dt
\end{equation}
with $X^{N,Z}_0=Z$. Here $(H^N_t,t\ge 0)$
is a $m$-dimensional Hermitian 
Brownian motion.  In other words, $H^N_t$ are a 
set of $m$ independent
matrix valued process whose matrix entries are given by
$$H^N_t(k,l)=\frac{B_{kl}(t)+i\tilde B_{kl}(t)}{\sqrt{2N}}
\, k<l, H^N_t(l,k)=\bar H^N_t(k,l), H^N_t(k,k)=\frac{B_{kk}(t)}{\sqrt{N}}$$
where the $B,\tilde B$ are independent standard Brownian motions.

A strong solution to \eqref{eqdiff}
 exists up to a possible time of
explosion (when $DV$ would stop being Lipschitz
eventually) since this equation can be seen
as a system of classical  stochastic
differential equation with $N(N+1)/2$ equations driven by 
independent Brownian motions 
with a polynomial drift.  Since this
drift derives from a strictly convex potential,
it is well known that the time of explosion is almost
surely infinite (which can also be show by the arguments
of the proof of Lemma \ref{exist})

Now, let us consider two solutions $X^{N,Z}$
and $X^{N,Z'}$ starting
from $Z$ and $Z'$ respectively.
Then, we have
$$d(X^{N,Z}_t-X^{N,Z'}_t)=-\half(DV(X^{N,Z}_t)-DV(X^{N,Z'}_t)) dt.$$
Hence, we can apply exactly the same arguments than in 
the proof of Theorem \ref{existtheo} to conclude
that 
$$\|X^{N,Z}_t-X^{N,Z'}_t\|_\infty\le e^{-ct}\|Z-Z'\|_\infty.$$
If we take $Z_N$ to be random with law  $\mu^N_V$
we get a stationary process so that
\begin{equation}
\label{cont}
\|X^{N,Z}_t-X^{N,0}_t\|_\infty\le e^{-ct}\|Z_N\|_\infty.
\end{equation}
Now, according to Brascamp-Lieb inequality (see 
its application on our particular case in \cite{GM2})
$$\mu^N_V(\|Z_N\|_\infty\ge x)\le e^{-a(c) N(x-x_0(c))}$$
with some $x_0(c)<\infty$ and $a_0(c)>0$.
Therefore, \eqref{cont} shows that for all $t\ge 0$,
\begin{equation}\label{contnorm}
\P(\|X^{N,0}_t\|_\infty\ge (e^{-ct} +1) K)
\le 2 e^{-a(c) N(K-x_0(c))}
\end{equation}
since both $Z_N$ and $X^{N,Z}_t$ have law $\mu_V^N$.

{\it Step 2: Uniform bounds on $X^{N,0}$.}
In this section we want to
show that we can also control uniformly the norm of $X^N_t=X^{N,0}_t$
for $t$ in a compact set. 
To do that let us remind that 
\begin{eqnarray*}
dX^N_t.X^N_t
&=&2X^N_t.dH^N_t-  DV(X^N_t).X^N_t dt +2dt\\
&\le& 2X^N_t.dH^N_t -c X^N_t.X^N_t dt +2dt-2DV(0).X^N_t dt\\
&\le& 2X^N_t.dH^N_t -\half c X^N_t.X^N_t dt +C
dt
\end{eqnarray*}
with $C=2+\frac{2}{c} DV(0).DV(0)$. 
Therefore, for any $p\ge 0$, since $X^N_t.X^N_t$
is a non negative matrix
for all $s\ge 0$, It\^o's calculus yield 
\begin{eqnarray}
\tr\left( (X^N_s.X^N_s)^p\right)&\le& 2p\int_0^s  \tr( (X^N_t.X^N_t)^{p-1}
(X^N_t.dH^N_t - \half c X^N_t.X^N_t dt +Cdt))\nonumber \\
&&+
2pN^{-1}
\sum_{k=0}^{p-1} \int_0^s \tr( (X^N_t.X^N_t)^k) \tr( (X^N_t.X^N_t)^{p-k-1}) dt
\nonumber\\
&&+
2pN^{-1}\sum_{k=0}^{p-2} \int_0^s \tr( (X^N_t.X^N_t)^k X^N_t ) 
\tr( (X^N_t.X^N_t)^{p-k-2} X^N_t ) dt\label{po}\\
\nonumber
\end{eqnarray}
Now, by Burkh\"older-Davis-Gundy inequalities and Chebychev's inequality,
there is a finite constant $\L$ such that 
for all $\e>0$,
\begin{multline*}
\P\left(\sup_{s\le T}
\left|\int_0^s  \tr( (X^N_t.X^N_t)^{p-1}
(X^N_t.dH^N_t)\right|\ge \e N^{\frac{1}{2}}\right)
\le \frac{\L}{\e^4 N^2}
E\left[\left(\frac{1}{N}\int_0^T \tr( (X^N_t.X^N_t)^{2p-1}
dt\right)^2\right]\\
\le\frac{\L T}{\e^4 N^2}\int_0^T 
E[\left(\frac{1}{N}\tr( (X^N_t.X^N_t)^{4p-2}
dt\right)] dt
\le \frac{\L T^2}{\e^4 N^2}( K^{8p} 
+ 4p \int_K^\infty x^{8p-4} e^{-a(c)N(x-x_0(c))} dx)\\
\le  \frac{\L T^2}{\e^4 N^2}( K^{8p} +\frac{(8p-4)!}{(a(c)N)^{8p-3}})
\le \frac{2\L T^2}{\e^4 N^2} K^{8p}
\end{multline*}
where  we used \eqref{contnorm}
and  chose $K=\max\{2x_0(c),1\}$, $p\le a(c)N K$.
Therefore, if we choose $\e$ so that $\e =(K+1)^{2p}T
$,
we see that if we set $$A(N,T)=\bigcap_{p\le a(c)N K} \left\{\sup_{s\le T}
\left|\int_0^s  \tr( (X^N_t.X^N_t)^{p-1}
(X^N_t.dH^N_t)\right|\le (K+1)^{2p} T  N^{\frac{1}{2}}\right\}$$
then by Borel Cantelli's  Lemma,
$$\limsup_{N,T\ra\infty} \P(A(N,T))=1.$$

Let now restrict ourselves to  the set
$\cap_{T\ge T_0}\cap_{N\ge N_0}A(N,T)$.
We let $$B(p,N,s):= N^{-1} \tr\left( (X^N_s.X^N_s)^p\right)$$
and observe that these non negative  real numbers  obey the relation
$$B(p,N,s)\le B(q,N,s)^{\frac{p}{q}},\qquad q\ge p.$$
We first control $B(p_0,N,s)$
by \eqref{po} which yields for $s\ge 0$, 
\begin{multline*}
B(p_0,N,s)\le 2p_0  (K+1)^{2p_0} N^{-\frac{1}{2}} s\vee T_0
+\int_0^s (-2cp B(p_0,N,t)+2 p_0(2p_0-1)
B(p_0,N,t)^{\frac{p_0-1}{p_0}}) dt \\
\le 2p_0  (K+1)^{2p_0}  N^{-\frac{1}{2}} s\vee T_0
+\int_0^s (-2cp B(p_0,N,t)+2 p_0(2p_0-1)
B(p_0,N,t)^{\frac{p_0-1}{p_0}}) dt \\
\leq 
 2p_0  (K+1)^{2{p_0}}  N^{-\frac{1}{2}}  s\vee T_0+ 
\int_0^s \Big\{ \big[ -2cp_0+a 2 p_0(2p_0-1)\big]
B(p_0,N,t) \\ +\left[
\Big(\frac{p_0-1}{p_0 a}\Big)^{p_0-1} -a \Big(\frac{p_0-1}{p_0 a}\Big)^{p_0}\right] \Big\}
dt\\ 
\end{multline*}
where we used that $x^{\frac{p_0-1}{p_0}}
\le a x+(\frac{p_0-1}{p_0 a})^{p_0-1} -a (\frac{p_0-1}{p_0 a})^{p_0}$
for all $a>0$ and $p_0\ge 1$.
Choosing $a$ so that $a (2 p_0+4 p_0(p_0-1))=cp_0$
we thus have found a finite constant $C(p_0,K)$
such that 
$$B(p_0,N,s)\le C(p_0,K)  s\vee T_0- cp_0 \int_0^s B(p_0,N,t) dt$$
which shows that $B(p_0,N,s)\le \max\{C(p_0,K)/cp_0, C(p_0,K) T_0\}  $
is uniformly bounded. 

We now bound $B(p,N,s)$ for $p\ge p_0$
and to this end  replace in \eqref{po} all
$B(q,N,s), q\le p_0$ by $ B(p_0,N,s)^{\frac{p_0}{q}}$.
We  show by induction over 
$p\ge p_0$ that $B(p,N,s)\le C_pC_0^{p}  $
with $C_0$ a finite constant depending on $p_0$ and $T_0$,
for all $s\ge 0$ and $p\le a(c)KN$.
Here $C_p$ denotes the Catalan numbers.
Indeed, this is satisfied for $p\le p_0$
and then \eqref{po} implies that
\begin{multline*}
B(p,N,s)\le 2p  (K+1)^{2p}  N^{-\frac{1}{2}}
s\vee T_0
+\int_0^s (-2cp B(p,N,t)+2 pB(p,N,t)^{\frac{p-1}{p}}) dt\\
+2p\sum_{k=0}^{p-1}\int_0^s C_{k} C_{p-k-1} C_0^{p-1} dt
\\ +2p\sum_{k=0}^{p-2}\int_0^s \Big[\left( C_{k+1} C_0^{k+1} \right)^{\frac
{k+\frac{1}{2}}{k+1}} \left( C_{p-k-1} C_0^{p-k-1} \right)^{\frac
{p-k-\frac{3}{2}}{p-k-1}}\Big]dt
\end{multline*}
where we used that $|\frac{1}{N}\tr( (X.X)^k X)|
\le \frac{1}{N}\tr((X.X)^{k+\frac{1}{2}})$.
Because $1\le C_k\le 4C_{k-1}$ and since $\sum_{k=0}^{p-1}
C_k C_{p-1-k}= C_{p}$ for all $k$ and $p$ we
conclude that
\begin{eqnarray*}
B(p,N,s)&\le& 2p  (K+1)^{2p}  N^{-\frac{1}{2}}
s\vee T_0
+\int_0^s (-cp B(p,N,t)+2 pB(p,N,t)^{\frac{p-1}{p}}) dt\\
&& +10 p C_0^{p-1} C_p s.\\
&\le& 2p  (K+1)^{2p}  N^{-\frac{1}{2}}
s\vee T_0
+\int_0^s \left(-\half cp B(p,N,t)+\left(\frac{2}{c}\right)^p \right) dt
 +10 p C_0^{p-1} C_p s.\\
\end{eqnarray*}
Thus, if $N$ is large enough so that $N^{-\frac{1}{2}}T_0\le 1$,
we get 
that
$$B(p,N,s)\le \frac{2 f(T_0)}{c}\left( (K+1)^{2p}  N^{-\frac{1}{2}}
+\left(\frac{2}{c}\right)^p+10  C_0^{p-1} C_p\right)$$
with $f(T_0)=T_0$ if $s\le T_0$ and
$f(T_0)=\frac{2}{c}$ if $s\ge T_0$. Note here that we used the
fact that we have a negative drift 
growing linearly with $p$
to cancel the multiplication by $p$.
We finally choose $C_0$ large
enough so that
$$\frac{2 f(T_0)}{c}\left( (K+1)^{2p}  N^{-\frac{1}{2}}
+\left(\frac{2}{c}\right)^p+10  C_0^{p-1} C_p\right)\le C_0^p$$
which we can always do.

Hence we have proved that
$$\sup_{t\ge 0}\|X^N_t.X^N_t\|_\infty\le
\min_{p\le a(c)KN} \sup_{t\ge 0}
(NB(p,N,s))^{1/ p} 
 \le 2C'(p_0)<\infty\quad \textrm{a.s.}$$
In other words, we have proved that $\| X^N_t\|_\infty$
is uniformly bounded almost  surely.

{\it Step 3: Convergence of the norm of $P(X^{N,0}_t)$ as $N$ goes
to infinity.}
To this end remark that since $X_t$ is always uniformly bounded by $M$
we can always assume $V$ is $\Ca^\infty$,
uniformly bounded 
and with uniformly Lipschitz cyclic gradient 
 (this amounts to change $V$
outside a place that the diffusion does not
see). 
We let  $\tilde V$ be  equal
to $V$ on operators with norm
bounded by $M$ and have 
uniformly Lipschitz  cyclic gradient.
 For instance, we take $\tilde V(X_1,\ldots, X_m)=
V( f(X_1),\ldots, f(X_m))$ with $f(x)= x$
on $|x|\le M$, $f(x)=x(1+(|x|-M)^4)^{-1}$ if $|x|>M$
(since the later is twice continuously differentiable with
uniformly bounded derivatives).

Now, by definition if 
we let
$$\phi_M(X,S)_t=S_t-\half \int_0^t D\tilde V(X_s)ds,$$ then $X_t$ can be 
expressed as an iterate
$$X_t=\phi_M(X,S)_t=\phi_M(.,S)^n(X)_t$$
for all integer numbers $n$ and $M$ greater than the uniform
norm on $X_t, t\ge 0$.
On the other hand, for two operator
valued processes $(X,Y)$
$$\|\phi_M(X,S)_t-\phi_M(Y,S)_t\|_\infty\le \half
\|D\tilde V\|_\La \int_0^t \|X_s-Y_s\|_\infty
ds$$
and so we get that
\begin{eqnarray*}
\|X_t-\phi(.,S)^n(S)_t\|_\infty&=&
\|\phi(X,S)_t-\phi(\phi^{n-1}(.,S)(S),S)_t\|_\infty\\
&\le&\|D\tilde V\|_\La\int_0^t\|X_s-\phi^{n-1}(S,.)(S)_s\|_\infty
ds\\
& \le &C\frac{\|D\tilde V\|_\La^n t^{n-1}}{2^n (n-1)!}\sup_{u\le t}\|X_u-S_u\|_\infty
\\
& \le &C\frac{\|D\tilde V\|_\La^n t^{n}}{2^n(n-1)!}\|D\tilde V\|_\infty
\end{eqnarray*}
for all $n\in\N$.

We next want to show that
the norm of $P(X^{N,0}_t)$ converges with
overwhelming
probability for any $P\in\cxm$.
To do that we approximate $X^{N,0}_t$
by
$\phi_M(.,H^N)^n(0)_t$ with $\phi$ as above

We claim that
$$\lim_{N\ra\infty}\lim_{M\ra\infty}\lim_{n\ra\infty}
\sup_{t\le T}\|\phi_M(.,H^N)^n(0)_t-X^N_t\|_\infty=0\qquad \mbox{a.s.}$$
Indeed, $\phi_M(.,H^N)$ is a contraction
for all $M$ finite and $(X^N_t, t\le T)$ is its unique fixed point
as long as $(X^N_t, t\le T)$ stays uniformly bounded 
by $M$. Since we have seen 
that almost surely $(X^N_t, t\le T)$ is  uniformly bounded,
the statement follows.

{\it Step 4: Convergence of the norm
of $\phi_M(.,H^N)$.}
For all $M$, $\phi_M(.,H^N)^n(0)_t$ can be approximated uniformly 
by a polynomial function
of $H^N$ on $\sup_{s\le t}\|H^N_s\|_\infty\le
L\}$ for some $L$ finite, which happens with 
 probability one
for some sufficiently large $L$.
We can thus use \cite{HT}
to conclude that the norm of
any polynomial in $\phi_M(.,H^N)^n(0)_t$
converges to its analog with $H^N$ replaced by $S$.

{\it Step 5: Conclusion.}
We have proved that for all $n$ and $M$
$$\lim_{N\ra\infty}
\|P(\phi_M(.,H^N)^n(0)_t)\|_\infty=\|P(\phi_M(.,S)^n(0)_t)\|_\infty\mbox{ a.s.}
$$
Since $X^N_t$ is uniformly bounded, for $N$ large enough, 
$$\lim_{M\ra\infty}\lim_{n\ra\infty}
\sup_{t\le T}\|\phi_M(.,H^N)^n(0)_t-X^N_t\|_\infty=0\qquad \mbox{a.s.}$$
implies 
$$\lim_{M\ra\infty}\lim_{n\ra\infty}
\sup_{t\le T}\|P(\phi_M(.,H^N)^n(0)_t)-P(X^N_t)\|_\infty=0\qquad \mbox{a.s.}$$
And finally we
have with overwhelming probability (where $Z_N$ with stationary law 
such that $\|Z_N\|\le K$)
 for $N$ large enough
$$\|X^{N,Z}_t-X^{N}_t\|_\infty\le e^{-2ct}\|Z_N\|_\infty\le e^{-ct} K.$$
Let $\e>0$ be fixed.
We fix $t$ so that $e^{-2ct} K=\e$.
With overwhelming probability, since $X^{N,Z}_t$
has the same law than $Z$, $\|X^{N,Z}_t\|_\infty\le K$
and so also $\|X^{N}_t\|_\infty\le \e+K$.
Hence, we have for any polynomial,
$$\lim_{N\ra\infty}
\|P(X^{N,Z}_t)-P(X^{N}_t)\|_\infty\le  C(K) \e  \quad\mbox{a.s.}$$
Now, we take
 $M,n$ large enough (greater than some
finite random integers) to assure that for $N$ large enough so that

$$\|P(\phi_M(.,H^N)^n(0)_t)-P(X^N_t)\|_\infty\le \e/3$$
and
finally
$$
|\|P(\phi_M(.,H^N)^n(0)_t)\|_\infty-\|P(\phi_M(.,S)^n(0)_t)\|_\infty|<\e/3.
$$
We have already seen that
$$\|P(\phi_M(.,S)^n(0)_t)-P(X_t)\|_\infty\le \e/3$$
Thus,  we
have proved
$$\big|
\|P(X^{N,Z}_t)\|_\infty-\|P(X^{Z}_t)\|_\infty\big|\le c'(K)\e.$$
This completes the proof since $X^{N,Z}_t$ (resp. $X^Z_t$) has the same law
that $Z_N$ (resp. $Z$).
\hfill$\square$

\end{document}